\documentclass[a4paper,11pt]{article}
\usepackage[utf8]{inputenc}
\usepackage{amsmath, amssymb}
\usepackage{geometry}
\usepackage{fancyhdr}
\usepackage{authblk}  
\usepackage{hyperref}
\usepackage{cite}
\usepackage{graphicx}

\usepackage{xcolor}
\usepackage{soul} 
\usepackage{booktabs}
\usepackage{subcaption}
\usepackage[title]{appendix}
\usepackage{siunitx}

\newcommand{\aver}[1]{\ensuremath{\{\!\{#1\}\!\}}}
\newcommand{\jump}[1]{\ensuremath{[\![ #1 ]\!]}}

\title{\textbf{Acoustic Propagation/Refraction Through Diffuse Interface Models}}

\author[1]{Abbas Ballout \thanks{Corresponding author: Abbas.Ballout@upm.es}}
\author[1]{Oscar A. Marino}
\author[1]{Gerasimos Ntoukas}
\author[1,2]{Gonzalo Rubio}
\author[1,2]{Esteban Ferrer}
\date{}

\affil[1]{ETSIAE-UPM-School of Aeronautics, Universidad Politécnica de Madrid, Plaza Cardenal Cisneros 3, E-28040 Madrid, Spain}
\affil[2]{Center for Computational Simulation, Universidad Politécnica de Madrid, Campus de Montegancedo, Boadilla del Monte, 28660 Madrid, Spain}

\begin{document}
\maketitle
\begin{abstract}
We present a novel approach for simulating acoustic (pressure) wave propagation across different media separated by a diffuse interface through the use of a weak compressibility formulation. Our method builds on our previous work on an entropy-stable discontinuous Galerkin spectral element method for the incompressible Navier-Stokes/Cahn-Hilliard system 
(Manzanero et al. (2020)), and incorporates a modified weak compressibility formulation that allows different sound speeds in each phase. We validate our method through numerical experiments, demonstrating spectral convergence for acoustic transmission and reflection coefficients in one dimension and for the angle defined by Snell's law in two dimensions. Special attention is given to quantifying the modeling errors introduced by the width of the diffuse interface. Our results show that the method successfully captures the behavior of acoustic waves across interfaces, allowing exponential convergence in transmitted waves. The transmitted angles in two dimensions are accurately captured for air-water conditions, up to the critical angle of $13^\circ$.  In a final example, we show a three-dimensional wave transmission from air into water to demonstrate the potential of this methodology for addressing general multiphase acoustic problems. This work represents a step forward in modeling acoustic propagation in incompressible multiphase systems, with potential applications to marine aeroacoustics.
\end{abstract}


\section{Introduction}
Modeling acoustic wave propagation through complex and heterogeneous media has widespread applications across multiple domains, including the medical field \cite{wang2022theoretical,doustikhah2024analysis}, aeroacoustics \cite{fatahian2020numerical,rutard2020large}, hydroacoustics \cite{portillo2024numerical,sezen2023marine,hynninen2023multiphase}, and in industry \cite{nair2024acoustic}. Accurately capturing the acoustic behavior in these scenarios often requires high-fidelity computational models that account for multiphase interactions and the presence of fluid interfaces.

Diffuse interface phase-field models are gaining traction in various multiphase modeling applications, offering an alternative to the well-established and mature sharp interface methods, namely volume of fluid (VOF) and level set (LS). This can be attributed to the fact that certain diffuse interface models can conserve mass while naturally incorporating surface tension effects. This allows them to maintain scalability on modern hardware without incurring additional algorithmic complications compared to VOF or LS \cite{mirjalili2017interface,li2023advances}. For instance, given the same required accuracy, it was demonstrated that a Cahn-Hilliard diffuse interface solver outperformed geometric VOF in terms of time-to-solution for canonical multiphase problems on second-order finite-difference structured grids \cite{mirjalili2019comparison}. Furthermore, numerical experiments using continuous Galerkin finite element discretization with adaptive mesh refinement have also compared LS to an Allen-Cahn diffuse interface model on a set of problems with large density ratios \cite{grave2022comparing}. The study concluded favorable results for LS in terms of accuracy, albeit the LS solver was more complicated to implement and incurred slightly additional computational overhead compared to the Allen-Cahn model. More interestingly, while diffuse interface models were originally developed to simulate phenomena in which the interface thickness is comparable to the physical length scales, they have recently shown a competitive advantage compared to VOF for scenarios where the physical length scales are much larger \cite{kuhl2022cahn}.  

Various methods exist for simulating acoustic propagation in multiphase flows, ranging from high-fidelity but computationally intensive models to comparatively more efficient yet simplified ones, all of which have been successfully applied to single-phase flows. The direct approach involves resolving the entire set of governing equations to obtain a high-fidelity pressure field, which is then used to compute the acoustic field (e.g., \cite{tajiri2010direct}). A less computationally demanding alternative is a two-step approach in which the flow field is first computed, and the linearized Euler equations are subsequently used to propagate pressure within the precomputed velocity field. This method has also been successfully applied to multiphase flows, as demonstrated in a series of works \cite{staab2015numerical,friedrich2018acoustics,friedrich2020towards}. However, since acoustic waves with small wave numbers require spatial discretization of comparable scale, the above methods become computationally prohibitive on large domains. In such scenarios, acoustic analogies offer a more feasible solution by propagating the pressure through a nonlinear equation based on surface pressure data from a high-fidelity solver. Among these models, the Ffowcs-Williams and Hawkings (FW-H) analogy is one of the most widely used, with successful applications in multiphase acoustic simulations reported in \cite{fatahian2020numerical,portillo2024numerical,hynninen2023multiphase,sezen2023marine,nair2024acoustic}.

In this work, we employ the high-order discontinuous Galerkin spectral element method (DGSEM) for spatial discretization \cite{kopriva2009implementing}. High-order methods are particularly appealing due to their low dissipative and dispersive errors \cite{manzanero2018dispersion,gassner2011comparison}. DGSEM, in particular, allows for the construction of provably stable discretization schemes, enabling kinetic energy-preserving and entropy-stable formulations for a range of equations, including the Euler equations \cite{gassner2016split}, compressible Navier-Stokes \cite{gassner2018br1}, the Spalart–Allmaras turbulence model for compressible Reynolds-Averaged Navier–Stokes equations \cite{lodares2022entropy}, the compressible multiphase Baer–Nunziato equations \cite{fraysse2016upwind, renac2019entropy, rai2021modelling}, and magnetohydrodynamics \cite{winters2016affordable,bohm2020entropy}. We have also successfully applied provably entropy-stable DGSEM schemes to the Cahn-Hilliard equation \cite{manzanero2020free}, incompressible Navier-Stokes with artificial compressibility \cite{manzanero2020entropyNS}, and the multiphase (two-phase and tri-phase) incompressible Navier-Stokes/Cahn-Hilliard (iNS/CH) system \cite{manzanero2020entropyNSCH,manzanero2022high} with p-adaptivity \cite{ntoukas2022entropy}.   
  
The predominant method for simulating incompressible flows has traditionally been splitting schemes \cite{guermond2006overview}. However, this approach yields a non-physical pressure field that merely enforces the incompressibility constraint, making it unsuitable for directly inferring the acoustic field. An alternative strategy is to employ an artificial or weak compressibility formulation \cite{chorin1997numerical, shen1997pseudo}. These formulations introduce an additional equation that couples velocity and pressure and have been successfully applied to multiphase flows \cite{manzanero2020entropyNSCH, shah2011numerical, kajzer2020weakly, ntouras2020coupled, nguyen2023review, orlando2024implicit, matsushita2021gas, yang2021weakly}. Initially, weak compressibility was used purely as a numerical technique to enforce incompressibility without a physical basis. However, with the appropriate coupling between the velocity and pressure fields, it can also be leveraged to propagate acoustic (pressure) waves at the correct speeds. By modifying the weak compressibility formulation, acoustic waves can be accurately modeled, as demonstrated in early single-phase flow applications dating back to 1993 \cite{manno1993developing}, and more recently in \cite{pont2018unified, jiang2022acoustic}. Nevertheless, to the best of the authors' knowledge, there are no reports in the literature on weak compressibility-enabled acoustic propagation in multiphase flows.

In this paper, we extend our previous work on DGSEM for the iNS/CH systems by introducing a novel approach for the direct acoustic wave propagation across different media using weak compressibility. This includes a modified weak compressibility formulation that allows for different sound speeds in each phase, separated by a diffuse interface. To ensure the physical soundness of our implementation, we validate its ability to replicate relevant physical phenomena, such as reflection, transmission, and refraction of acoustic waves according to Snell's law. We also examine the DGSEM's spectral convergence and the modeling errors introduced by the diffuse interface.

The rest of the paper is organized as follows. Section \ref{sec:Continuous-setting} presents the governing equations in the continuous setting, followed by the DGSEM machinery employed for spatial discretization in general, and the numerical fluxes responsible for the propagation of the speed of sound in particular in section \ref{sec:Discretization}. Section \ref{sec:Experiments} describes our numerical experiments, first investigating one-dimensional acoustic wave transmission and reflection, then two-dimensional wave propagation with Snell's law validation, and finally, we extend to three-dimensional propagation of spherical waves across a flat interface, assuming air water conditions. Section \ref{sec:Conclusion} concludes the study and examines the implications of the proposed approach for multiphase acoustic modeling.

\section{Continuous setting}\label{sec:Continuous-setting}
This section summarizes the mathematical model used to simulate acoustic wave propagation within two media. We extend our previous work of a two-phase entropy-stable DGSEM iNS/CH model with artificial compressibility, valid for high density ratios \cite{manzanero2020entropyNSCH}, to account for different sound speeds in each phase by modifying the artificial compressibility equation. 
Throughout this work, scalar variables are explicitly denoted by lowercase characters, and vectors are denoted by bold lowercase characters. 

The concentration of each phase $c \in [0,1]$ is governed by the diffuse Cahn-Hilliard equation with advection \cite{cahn1958free,cahn1959free, berthier2001phase}:
\begin{equation}
    \partial_t c + \nabla \cdot \left(c \boldsymbol{u} \right) = M_0 \nabla^2 \mu ,
    \label{eq:Cahn-Hilliard}
\end{equation}
 with $\mu$ as the chemical potential, given by:
\begin{equation}
      \mu = \dfrac{d f_0}{dc} -\dfrac{3}{2}\sigma \varepsilon \nabla^2 c,  \quad f_0 = \dfrac{12 \sigma}{\varepsilon} c^2(1-c)^2,
    \label{eq:chemical-potential}
\end{equation}
such that the surface tension $\sigma$, the diffuse interface width $\varepsilon$, and an additional chemical characteristic time $t_{CH}$ parameter relate to the mobility parameter $M_0$ through the relation:
\begin{equation}
    M_0 = \dfrac{\varepsilon^2}{\sigma t_{CH}}.
    \label{eq:mobility}
\end{equation}
The advective velocity $\boldsymbol{u}=(u,v,w)$ in equation \eqref{eq:Cahn-Hilliard} couples the Cahn-Hilliard equation to the Navier-Stokes equation, given by: 
\begin{equation}
    \sqrt{\rho} \partial_t \left( \sqrt{\rho} \boldsymbol{u} \right)
    +\nabla \cdot \left( \dfrac{1}{2} \rho \boldsymbol{u} \boldsymbol{u}  \right) 
    +\dfrac{1}{2} \rho \boldsymbol{u} \cdot \nabla \boldsymbol{u} 
    + c \nabla \mu
    = - \nabla p  
    + \nabla \cdot \left(\eta  (\nabla \boldsymbol{u} + \nabla \boldsymbol{u}^T) \right) 
    + \rho \boldsymbol{g} ,
    \label{eq:Navier-Stokes}
\end{equation}
where $p$ is pressure, $\boldsymbol{g}$ is gravity, $\rho$ is density, and $\eta$ is viscosity. While we have previously used weak compressibility to impose the incompressibility constraint, we now utilize it to propagate the pressure at the correct acoustic speed. Assuming isentropic behavior and that the acoustic speed is much larger than the advective velocity (low Mach numbers), we write the pressure equation as \cite{clausen2013entropically}:  
\begin{equation}
    \partial_t p + \rho c_s^2 \nabla \cdot \boldsymbol{u} = 0 ,
    \label{eq:Pressure}
\end{equation}
with $c_s$ being the speed of sound. Within the diffuse interface, density $\rho$, viscosity $\eta$, and the speed of sound $c_s$ are interpolated according to the concentration field as:
\begin{equation}
   (.) = (.)_1 c + (.)_2 (1-c),
\end{equation}
where $(.)$ represents a general property, and subscripts $1$ and $2$ correspond respectively to each phase. Equations \eqref{eq:Cahn-Hilliard}, \eqref{eq:Navier-Stokes}, and \eqref{eq:Pressure} form the multiphase iNS/CH system, which will be shown to accurately model acoustic wave propagation in two phases.  

\section{Discretization}\label{sec:Discretization}

In this section, we present the development of the DGSEM discretization scheme, which will be used for the spatial semi-discretization of the iNS/CH system introduced in the previous section. We also present the non-conservative fluxes that involve the speed of sound. As this section is not intended to provide an exhaustive introduction to DGSEM, readers are encouraged to consult the references cited in the introduction.

\subsection{Curvilinear mapping of operators}
\label{subsec:curvelinear-mapping}
Consider a physical domain $\Omega \subseteq \mathbb{R}^3$ partitioned into $K$ non-overlapping hexahedral elements, denoted by $e_k$ for $k = 1, 2, \ldots, K$. Each element $e_k$ is mapped from a reference element $E = [-1, 1]^3$ using a transfinite mapping $\boldsymbol{X}^e$. This mapping relates the physical coordinates $\boldsymbol{x} = (x^1, x^2, x^3)$ to the reference coordinates $\boldsymbol{\xi} = (\xi^1, \xi^2, \xi^3)$ through:
\begin{equation}
\boldsymbol{\chi} = \boldsymbol{X}^e(\boldsymbol{\xi}).
\end{equation}
For simplicity, the superscript $e$ will be omitted in the subsequent expressions. From the transformation $\boldsymbol{\chi} = \boldsymbol{X}(\boldsymbol{\xi})$, the covariant basis vectors are defined as:
\begin{equation}
\boldsymbol{a}_i = \frac{\partial \boldsymbol{X}}{\partial \xi^i}, \quad i = 1, 2, 3,
\end{equation}
and the contravariant basis vectors $\boldsymbol{a}^i$ are determined by:
\begin{equation}
\boldsymbol{a}^i = \nabla \xi^i = \dfrac{1}{J}(\boldsymbol{a}_j \times \boldsymbol{a}_k  ), \quad (i,j,k) \,\, \text{cyclic}, 
\label{eq:contravariat-bases}
\end{equation}
where $J$ is the Jacobian determinant of the mapping given by:
\begin{equation}
J = \boldsymbol{a}_1 \cdot (\boldsymbol{a}_2 \times \boldsymbol{a}_3).
\end{equation}
The operators are transformed from real space into the reference coordinates via the metric matrix $\boldsymbol{M}$, so that the divergence and gradient operators become \cite{gassner2018br1}:
\begin{align}
\nabla\cdot \boldsymbol{f} = \frac{1}{J} \nabla_\xi \cdot \left(\boldsymbol{M}^T\boldsymbol{f}\right), \quad
\nabla {f} = \frac{1}{J}  \boldsymbol{M} \nabla_\xi {f}, \quad
\nabla \boldsymbol{f} = \frac{1}{J} \boldsymbol{M}  \nabla_\xi \boldsymbol{f}, 
\label{eq:Operators}
\end{align}
where $\boldsymbol{M} = (J \boldsymbol{a}^1, J \boldsymbol{a}^2, J \boldsymbol{a}^3 )$. The extension to block operators on block entries of size $n$, is then readily available by extending the definition of the metric matrix instead to be: 
\begin{equation}
\mathcal M_n = \left(\begin{array}{ccc}
Ja^{1}_{1}\boldsymbol{I}_{n} & Ja^{2}_{1}\boldsymbol{I}_{n} & Ja^{3}_{1}\boldsymbol{I}_{n}\\
Ja^{1}_{2}\boldsymbol{I}_{n} & Ja^{2}_{2}\boldsymbol{I}_{n} & Ja^{3}_{2}\boldsymbol{I}_{n}\\
Ja^{1}_{3}\boldsymbol{I}_{n} & Ja^{2}_{3}\boldsymbol{I}_{n} & Ja^{3}_{3}\boldsymbol{I}_{n}\end{array}\right),
\label{eq:metrics:metrics-matrix}
\end{equation}
where $\boldsymbol{I}_n$ is an $n \times n$ identity matrix. The operator transformations in \eqref{eq:Operators} remain valid under this extension. This block-type formulation proves useful for handling systems of equations, as in the iNS/CH system. 

\subsection{Polynomial approximation and DGSEM formulation}
\label{subsec:DGSEM}
Within each element, the unknown vector variables $\boldsymbol{q}$ are approximated as an order $N$ polynomial, which is allowed to vary from element to element:
\begin{equation}
  \boldsymbol{q}(\boldsymbol{{\xi}},t) \approx \Tilde{\boldsymbol{q}} (\boldsymbol{\xi},t) = \sum_{i,j,k=0}^{N} \Tilde{\boldsymbol{q}}_{ijk}(t) \, l_i(\xi^1) \, l_j(\xi^2) \, l_k(\xi^3),
  \label{eq:interp}
\end{equation}
where the Lagrange polynomials:
\begin{equation}
l_i(\xi) = \prod_{\substack{j=0\\j\neq i}}^{N}\frac{\xi-\xi_j}{\xi_i-\xi_j},
\end{equation}
are defined at Gauss–Lobatto nodes $\{\xi_i\}_{i=0}^{N}$. The choice of nodes is not limited to those of Gauss-Lobatto; however, Gauss-Lobatto nodes satisfy the SBP-SAT property \cite{carpenter2013high}, which allows us to discretely mirror the continuous stability framework and establish a discrete entropy law, as shown in our earlier work on the iNS/CH system \cite{manzanero2020entropyNSCH,ntoukas2022entropy}.

Both geometry and metric terms are interpolated in $\mathbb{P}^N$, 
and the metric identities for the contravariant coordinate vectors \cite{kopriva2006metric}:
\begin{equation}
  \sum_{i=1}^{3}\frac{\partial (J\, a_n^i)}{\partial \xi^i} = 0,\quad n=1,2,3,
\end{equation}
which ensures that the discrete free-stream preservation holds for the continuous setting. However, they do not hold discretely if the contravariant vectors are constructed according to \eqref{eq:contravariat-bases}, since the product $a_j\times a_k$ is a polynomial of order $2N$. Therefore, the discrete contravariant bases are constructed through \cite{kopriva2006metric}:
\begin{equation}
 J a_{n}^{i} = -\hat{x}_i \cdot\nabla_{\xi}\times\left(\mathbb I^{N}\left(\mathcal X_{l}\nabla_{\xi}\mathcal X_{m}\right)\right)\in\mathbb{P}^{N},~~~i=1,2,3,~~n=1,2,3,~~(n,m,l)\text{ cyclic},
\label{eq:dg:contravariant-basis}
\end{equation}
where $\mathbb{I}^N$ denotes the polynomial interpolation operator. 

Volume integrals are approximated using tensor-product Gauss-Lobatto quadratures defined on the same Gauss-Lobatto nodes:
\begin{equation}
\langle f,g \rangle_{E,N} \approx \sum_{m,n,l=0}^{N} w_m\,w_n\,w_l \, f_{mnl}\, g_{mnl},
\end{equation}
which defines an inner product in the reference element $E$. The same applies to surface integrals.

\subsection{Spatial discretization}
\label{subsec:spatial-discretization}
With the DGSEM discretization framework established, we can now construct the discrete equivalent of our continuous formulation. Following the notation in \cite{gassner2018br1}, we begin by expressing equations \eqref{eq:Cahn-Hilliard}, \eqref{eq:Navier-Stokes}, and \eqref{eq:Pressure} more compactly: 
\begin{equation}
    \underline{\boldsymbol{M}} \partial_t \boldsymbol{q} + \nabla \cdot \boldsymbol{F}_e (\boldsymbol{q}) + \sum_{m=1}^{5} \boldsymbol{\Phi}_m (\boldsymbol{q}) \cdot \nabla w_m = \nabla \cdot \boldsymbol{F}_{\nu} (\nabla \boldsymbol{w}) + \boldsymbol{s}(\boldsymbol{q})
    \label{eq:compactiNS/CH-system},
\end{equation}
with state vector $\boldsymbol{q}=(c,\sqrt{\rho}\boldsymbol{u},p)$, gradient variables vector $\boldsymbol{w}= ( \mu, u, v,w, p)$, mass matrix:
\begin{equation}
    \underline{\boldsymbol{M}} =  
    \begin{pmatrix}
  1 & \boldsymbol{0} & 0  \\
  \boldsymbol{0} & \sqrt{\rho} \boldsymbol{I}_3 & \boldsymbol{0} \\
  0 & \boldsymbol{0} & 1\\
\end{pmatrix},
\end{equation}
inviscid fluxes $\boldsymbol{F}_e (\boldsymbol{q})=(\boldsymbol{f}_{e,1},\boldsymbol{f}_{e,2}, \boldsymbol{f}_{e,3} )$:
\begin{equation}
    \boldsymbol{f}_{e,1} = \boldsymbol{f}_{e} =  
    \begin{pmatrix}
        cu \\
        \frac{1}{2} \rho u^2 + p \\
         \frac{1}{2} \rho u v \\
         \frac{1}{2}\rho u w \\
        0
    \end{pmatrix}, 
     ~   \boldsymbol{f}_{e,2} = \boldsymbol{g}_{e} =  
    \begin{pmatrix}
        cv \\
         \frac{1}{2} \rho u v  \\
         \frac{1}{2} \rho u v^2 + p \\
         \frac{1}{2} \rho v w \\
        0
    \end{pmatrix}, 
     ~   \boldsymbol{f}_{e,3} = \boldsymbol{h}_{e} =  
    \begin{pmatrix}
        cw \\
         \frac{1}{2} \rho u w  \\
         \frac{1}{2} \rho v w  \\
         \frac{1}{2} \rho  w^2 + p \\
        0
    \end{pmatrix}
\end{equation}
non-conservative term coefficients:$\boldsymbol{\Phi}_m(\boldsymbol{q})$:
\begin{equation}
    \boldsymbol{\Phi}_1 =
        \begin{pmatrix}
        \boldsymbol{0} \\
        c \boldsymbol{e}_1  \\
        c \boldsymbol{e}_2  \\
        c \boldsymbol{e}_3 \\
        \boldsymbol{0}
    \end{pmatrix}, ~
        \boldsymbol{\Phi}_2 = 
    \begin{pmatrix}
        \boldsymbol{0} \\
        \frac{1}{2} \rho \boldsymbol{u}  \\
        \boldsymbol{0}  \\
        \boldsymbol{0} \\
        \rho c_s^2 \boldsymbol{e}_1
    \end{pmatrix}, ~
        \boldsymbol{\Phi}_3 = 
    \begin{pmatrix}
        \boldsymbol{0} \\
        \boldsymbol{0}  \\
        \frac{1}{2} \rho \boldsymbol{u}  \\
        \boldsymbol{0} \\
        \rho c_s^2 \boldsymbol{e}_2
    \end{pmatrix}, ~
        \boldsymbol{\Phi}_4 = 
    \begin{pmatrix}
        \boldsymbol{0} \\
        \boldsymbol{0}  \\
        \boldsymbol{0}  \\
        \frac{1}{2} \rho \boldsymbol{u} \\
        \rho c_s^2 \boldsymbol{e}_3
    \end{pmatrix}, ~ 
        \boldsymbol{\Phi}_5 = 
    \begin{pmatrix}
        \boldsymbol{0} \\
        \boldsymbol{0}  \\
        \boldsymbol{0}  \\
        \boldsymbol{0} \\
        \boldsymbol{0}
    \end{pmatrix}, ~ 
\end{equation}
viscous fluxes: $\boldsymbol{F}_{\nu} (\nabla \boldsymbol{w})=
(\boldsymbol{f}_{\nu,1}, \boldsymbol{f}_{\nu,2},\boldsymbol{f}_{\nu,3})$,
\begin{equation}
    \boldsymbol{f}_{\nu,1} = \boldsymbol{f}_{\nu} =  
    \begin{pmatrix}
        M_0 \partial_x \mu \\
        2 \eta \boldsymbol{S}_{11} \\
        2 \eta \boldsymbol{S}_{21} \\
        2 \eta \boldsymbol{S}_{31} \\
        0
    \end{pmatrix}, 
     ~
    \boldsymbol{f}_{\nu,2} = \boldsymbol{g}_{\nu} =  
    \begin{pmatrix}
        M_0 \partial_y \mu \\
        2 \eta \boldsymbol{S}_{12} \\
        2 \eta \boldsymbol{S}_{22} \\
        2 \eta \boldsymbol{S}_{32} \\
        0
    \end{pmatrix}, 
     ~
    \boldsymbol{f}_{\nu,3} = \boldsymbol{h}_{\nu} =  
    \begin{pmatrix}
        M_0 \partial_z \mu \\
        2 \eta \boldsymbol{S}_{13} \\
        2 \eta \boldsymbol{S}_{23} \\
        2 \eta \boldsymbol{S}_{33} \\
        0
    \end{pmatrix}, 
\end{equation}
and source terms $\boldsymbol{s}(\boldsymbol{q})=(0,\rho \boldsymbol{g},0)$, where $\boldsymbol{S}=\frac{1}{2}\left( \nabla \boldsymbol{u} + \nabla \boldsymbol{u}^T \right)$. We rewrite equation \eqref{eq:compactiNS/CH-system} given by a fourth-order operator, as a system of four first-order equations by introducing the auxiliary variables $\boldsymbol{G} = \nabla \boldsymbol{w}$ and $g = \nabla c$ so that:
\begin{equation}
\begin{aligned}
    \underline{\boldsymbol{M}} \partial_t \boldsymbol{q} &+ \nabla \cdot \boldsymbol{F}_e (\boldsymbol{q}) + \sum_{m=1}^{5} \boldsymbol{\Phi}_m (\boldsymbol{q}) \cdot \nabla w_m = \nabla \cdot \boldsymbol{F}_{\nu} (\boldsymbol{G}) + \boldsymbol{s}(\boldsymbol{q}), \\
    \boldsymbol{G} &= \nabla \boldsymbol{w}, \\
    \mu &= \dfrac{d f_0}{dc} -\dfrac{3}{2}\sigma \varepsilon \nabla \cdot \boldsymbol{g},  \\
    \boldsymbol{g} &= \nabla c.
    \label{eq:1storderNS/CH-systemo}
\end{aligned}
\end{equation}
By transforming the above set of equations to the reference space using \eqref{eq:Operators}, multiplying with test functions $\boldsymbol{\phi}_q$, $\boldsymbol{\Phi}_g$, $\phi_{\mu}$ and $\boldsymbol{\phi}_c$, inserting the discontinuous polynomial approximation \eqref{eq:interp}, and integrating by parts, we arrive at the discrete weak form:
\begin{equation}
\begin{split}
 \langle J  \underline{\boldsymbol{M}} \partial_t \tilde{\boldsymbol{q}} , \tilde{\boldsymbol{\phi}}_q \rangle_{E,N} 
&+\int\limits_{\partial E, N} \tilde{\boldsymbol{\phi}}_q^T \left( \mathcal{M}_5^T \tilde{\boldsymbol{F}}_e^* + \sum_{m=1}^{5} \left( \tilde{\boldsymbol{\Phi}}_m \mathcal{M}_5 \tilde{w}_m \right)^{\diamond} - \mathcal{M}_5^T \tilde{\boldsymbol{F}}_\nu^* \right) \cdot \hat{\boldsymbol{n}} dS_{\xi} 
-  \langle \mathcal{M}_5^T \tilde{\boldsymbol{F}}_e , \nabla_\xi \tilde{\boldsymbol{\phi}}_q   \rangle_{E,N} \\
&- \sum_{m=1}^{5} \langle \mathcal{M}_5 \tilde{w}_m, \nabla_\xi \cdot (\tilde{\boldsymbol{\phi}}_q^T \tilde{\boldsymbol{\Phi}}_m ) \rangle_{E,N}
=- \langle \mathcal{M}_5^T \tilde{\boldsymbol{F}}_\nu , \nabla_\xi \tilde{\boldsymbol{\phi}}_q   \rangle_{E,N} 
+ \langle J \tilde{\boldsymbol{s}}, \tilde{\boldsymbol{\phi}}_q  \rangle_{E,N} \\
\langle J \tilde{\boldsymbol{G}}, \tilde{\boldsymbol{\Phi}}_g  \rangle_{E,N}
&=\int\limits_{\partial E, N} \mathcal{M}_5 \tilde{\boldsymbol{w}}^{*,T} \tilde{\boldsymbol{\Phi}}_g \cdot  \hat{\boldsymbol{n}} dS_{\xi}
- \langle \mathcal{M}_5 \tilde{\boldsymbol{w}}, \nabla_{\xi} \cdot \tilde{\boldsymbol{\Phi}}_g  \rangle_{E,N} \\
\langle J \tilde{\mu}, \tilde{\phi}_\mu  \rangle_{E,N}
&= \langle J \dfrac{d \tilde{f}_0}{dc}, \tilde{\phi}_\mu  \rangle_{E,N}
- \int\limits_{\partial E, N}  \dfrac{3}{2}\sigma \varepsilon \tilde{\phi}_\mu \mathcal{M}_5 \tilde{\boldsymbol{g}}^* \cdot  \hat{\boldsymbol{n}} dS_{\xi} 
+ \langle \dfrac{3}{2}\sigma \varepsilon \mathcal{M}_5 \tilde{\boldsymbol{g}}, \nabla_{\xi} \tilde{\phi}_\mu \rangle_{E,N} \\
\langle J \tilde{\boldsymbol{g}}, \tilde{\boldsymbol{\phi}}_c  \rangle_{E,N} 
&= \int\limits_{\partial E, N}  M \tilde{c}^* \tilde{\boldsymbol{\phi}}_c  \cdot  \hat{\boldsymbol{n}} dS_{\xi} 
- \langle \tilde{c}, \nabla_\epsilon \cdot M \tilde{\boldsymbol{\phi}}_c  \rangle_{E,N}
\end{split}
\label{eq:discrete-weak-form}
\end{equation}
where $\hat{\boldsymbol{n}}$ is the unit normal, and  $dS_{\xi}$ is the surface differential on each face of the reference element $E$. 
The semi-discrete formulation is completed by introducing numerical fluxes to enforce inter-element coupling. The star superscripts refer to numerical fluxes for the conservative inviscid and conservative viscous fluxes, whose expression yields an entropy stable scheme, and are given in our previous work \cite{manzanero2020entropyNSCH}. The non-conservative diamond fluxes, which incorporate the spatially varied speed of sound, will be discussed in the next section. 

\subsection{Non-conservative diamond fluxes}
\label{subsec:non-conserv-fluxes}
In what follows, the tilde notation will be dropped, and all variables are understood to be in the discrete setting. We first rewrite the surface integral of the non-conservative inviscid flux in real space:
\begin{equation}
    \int\limits_{\partial E, N} {\boldsymbol{\phi}}_q^T \sum_{m=1}^{5} \left( {\boldsymbol{\Phi}}_m \mathcal{M}_5 {w}_m \right)^{\diamond} \cdot \hat{\boldsymbol{n}} dS_{\xi} =
    \int\limits_{\partial e, N} {\boldsymbol{\phi}}_q^T \sum_{m=1}^{5} \left( {\boldsymbol{\Phi}}_m {w}_m \right)^{\diamond} \cdot \boldsymbol{n} dS,
\end{equation}
where $\boldsymbol{n} = (n_x,n_y,n_z)$ is the unit normal, and $dS$ is the surface differential on each face of the real element $e$. The rotational invariance of the flux \cite{manzanero2020entropyNSCH, toro2013riemann} allows us to express the non-conservative terms as: 
\begin{equation}
    \sum_{m=1}^{5} \left( {\boldsymbol{\Phi}}_m {w}_m \right) \cdot \boldsymbol{n} = \boldsymbol{T}^T 
      \begin{pmatrix}
      0 \\
      \frac{1}{2} \rho u_n^2 + \mu c \\
      \frac{1}{2} \rho u_n v_{t1} \\
      \frac{1}{2} \rho u_n v_{t2} \\
      \rho c_s^2 u_n
      \end{pmatrix},~
      \boldsymbol{T}  = 
      \begin{pmatrix}
      1 & 0         & 0   & 0   & 0 \\
      0 & n_x      & n_y & n_z & 0 \\
      0 & t_{1,x} & t_{1,y} & t_{1,z} & 0 \\
      0 & t_{2,x} & t_{2,y} & t_{2,z} & 0 \\
      0 & 0 & 0 & 0 & 1
      \end{pmatrix},
\end{equation}
such that $\boldsymbol{T}$ is the rotation matrix, and $\boldsymbol{t}_1$, $\boldsymbol{t}_2$ are the two face tangent unit vectors. Multiplying the state vector \(\boldsymbol{q}\) by the rotation matrix $\boldsymbol{T}$ results in the face normal state vector $\boldsymbol{q}_n$:
\begin{equation}
    \boldsymbol{q_n} = \boldsymbol{T} \boldsymbol{q} = (c, \sqrt{\rho} u_n, 
    \sqrt{\rho} v_{t1}, \sqrt{\rho} v_{t2}, p)
\end{equation}
where $u_n= \boldsymbol{u} \cdot \boldsymbol{n}$ is the normal velocity, and $v_{ti}= \boldsymbol{u}\cdot \boldsymbol{t}_i ~ (i=1,2) $ are the two tangential velocities. 
Two options for the non-conservative diamond fluxes are provided:
\begin{enumerate}
    \item \emph{Entropy--conserving central fluxes:} The numerical flux is adapted from \cite{bohm2020entropy} and written as:
    \begin{equation}
    \sum_{m=1}^{5} \left( {\boldsymbol{\Phi}}_m {w}_m \right)^{\diamond} \cdot \boldsymbol{n} = \boldsymbol{T}^T
      \begin{pmatrix}
      0 \\
      \frac{1}{2} \rho u_n \aver{u_n} + c \aver{\mu} \\
      \frac{1}{2} \rho u_n \aver{v_{t1}} \\
      \frac{1}{2} \rho u_n \aver{v_{t2}} \\
      \rho c_s^2 \aver{u_n}
      \end{pmatrix},
      \label{eq:central-flux}
    \end{equation}
    where $\aver{}$ denotes the average operator for any quantity on the left and on the right, $\aver{\cdot} = \frac{1}{2} \left (\cdot_L + \cdot_R \right) $
    \item \emph{Entropy--stable Exact Riemann Solver (ERS):} The star region solution is computed from the left/right states via \cite{bassi2018artificial}:
    \begin{equation}
    \begin{split}
      u_n^\star &= \frac{p_L-p_R+\rho_L u_{nL}\lambda_L^+ -\rho_R u_{nR}\lambda_R^-}{\rho_L \lambda_L^+ - \rho_R \lambda_R^-}, \quad
      p^\star = p_L + \rho_L \lambda_L^+ \left( u_{nL}-u_n^\star\right), \\[1em]
    \rho^{\star} &= \left\{
    \begin{array}{ccc}
    \rho_{L}^{\star} & \text{ if } & u_n^\star \geqslant 0 \\
    \rho_{R}^{\star} & \text{ if } & u_n^\star < 0 
    \end{array}\right., \quad
    \rho_L^{\star} = \frac{\rho_L \lambda_L^{+}}{u_n^{\star}-\lambda_L^-}, \quad \rho_R^{\star} = \frac{\rho_R \lambda_R^{-}}{u_n^{\star}-\lambda_R^+}, \quad
    v_{ti}^{\star} = \left\{
    \begin{array}{ccc}
    v_{tiL} & \text{ if } & u_n^\star \geqslant 0 \\
    v_{tiR} & \text{ if } & u_n^\star < 0 
    \end{array}\right. ,
    \end{split}
    \end{equation}
    with eigenvalues:
    \begin{equation}
    \lambda_L^{\pm} = \frac{u_{nL} \pm a_L}{2}, \quad \lambda_R^{\pm} = \frac{u_{nR} \pm a_R}{2}, \quad a = \sqrt{u_n^2 + 4 c_s^2}.
    \end{equation}
    For non-conservative terms, we choose the diamond fluxes, 
\begin{equation}
    \begin{split}
    \sum_{m=1}^{5} \left( {\boldsymbol{\Phi}}_m {w}_m \right)^{\diamond} \cdot \boldsymbol{n} &= 
    {\boldsymbol{\Phi}}_1 \aver{{w}_1} \cdot \boldsymbol{n} + \sum_{m=2}^{5}\left({\boldsymbol{\Phi}}_m^* {w}_m^* +
    {\boldsymbol{\Phi}}_m {w}_m - 
    {\boldsymbol{\Phi}}_m^* {w}_m \right)\cdot \boldsymbol{n}\\
    &= \boldsymbol{T}^{T}
    \begin{pmatrix}
    0 \\
    \frac{1}{2}\rho^{*}u_{n}^{*,2} + \frac{1}{2}\rho u_{n}^{2} - 
    \frac{1}{2}\rho^{*}u_n^{*}u_n + c\aver{\mu} \\
    \frac{1}{2}\rho^{*}u_{n}^{*}v_{t1}^{*} + \frac{1}{2}\rho u_{n}v_{t1} - 
    \frac{1}{2}\rho^{*}u_n^{*}v_{t1} \\
    \frac{1}{2}\rho^{*}u_{n}^{*}v_{t2}^{*} + \frac{1}{2}\rho u_{n}v_{t2} - 
    \frac{1}{2}\rho^{*}u_n^{*}v_{t2} \\
    \rho c_s^2 u_n^{\star}
    \end{pmatrix}.
    \end{split}
    \label{eq:ERS-flux}
\end{equation}
\end{enumerate}

Both the numerical central flux (Equation \eqref{eq:central-flux}) and the numerical ERS flux (Equation \eqref{eq:ERS-flux}) were shown to be entropy stable on both interior and boundary facets, as demonstrated in detail \cite{manzanero2020entropyNSCH}. We confirm that this remains the case as reported in \ref{sec:Appendix-B}.
In the next section, we evaluate the model’s capability to reproduce acoustic phenomena through a series of numerical experiments. These tests are designed to rigorously assess the method's performance in capturing key acoustic features, such as wave propagation, reflection, transmission, and refraction, by comparing the computed results against theoretical predictions in one, two, and three dimensions. Additionally, since all validations are based on analytical solutions that exclude source terms, the gravity term is therefore neglected.

\section{Numerical experiments}\label{sec:Experiments}
All simulations use an explicit third-order low-storage Runge-Kutta time-stepping scheme. 
The entropy-stable ERS numerical flux is utilized for both inviscid conservative and non-conservative fluxes, while BR1 is applied for the viscous flux, as they have been demonstrated to ensure entropy stability on Gauss-Lobatto nodes \cite{manzanero2020entropyNSCH}. The implementation was developed within and uses HORSES3D \cite{ferrer2023high}, a high-order DGSEM solver. While HORSES3D allows for anisotropic polynomial orders in each element and across the domain, here we use a constant polynomial order uniformly. 

\subsection{One dimension - Reflection and transmission}
\label{subsec:1D}

\begin{figure}[h] 
    \centering
    \includegraphics[width=0.75\textwidth]{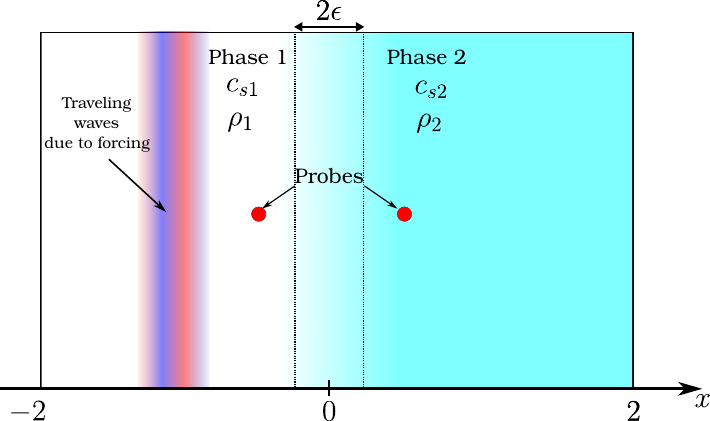}
    \caption{Problem setup for the 1D problem (not to scale).}
    \label{fig:1D-setup}
\end{figure}

This numerical experiment examines the transmission and reflection of an acoustic plane wave across two distinct phases. It also aims to quantify the modeling errors introduced by the diffuse interface. To this end, we present a one-dimensional numerical setup shown in Fig.~\ref{fig:1D-setup}, consisting of two media separated by a diffuse interface of thickness $2 \varepsilon$, within a domain that spans $\SI{4}{\meter}$ along the x-axis. The initial condition for the concentration field, which centers the interface at the origin, is defined as follows:
\begin{equation}
    c_0 = 1.0 - 0.5\left(1+tanh \left( \dfrac{2 x} {\varepsilon} \right)\right).
\end{equation}
Waves are excited within the domain by adding a forcing term to the right-hand side of the pressure equation \eqref{eq:Pressure}. The forcing term represents an oscillating Gaussian pulse with frequency $f$, width $b$, centered around $x_0$:
\begin{equation}
    s = cos( 2\pi f  t) exp \left(\left(\dfrac{x-x_0}{b} \right)^2\right) ;
    \label{eq:source}
\end{equation}
$b=\SI{0.01}{\meter}$ and $x_0 = \SI{-0.55}{\meter}$ for the remainder of the manuscript unless otherwise mentioned. The left and right walls are modeled as impermeable, no-slip boundaries. However, none of the simulations ran for a long enough duration for the signal to reach either boundary within the time required to complete the simulation. For the concentration field and chemical potential, no-flux boundary conditions are applied. Specifically, the boundary conditions are given by:
\begin{equation}
    \boldsymbol{n} \cdot \nabla c = 0 \quad \text{and} \quad \boldsymbol{n} \cdot \nabla \mu = 0.
    \label{eq:no-flux}
\end{equation}

Two probes are placed on both sides of the interface to monitor the pressure signals. The probe on the right measures the transmitted pressure wave $ p_t $, while the probe on the left records the sum of the incident and reflected pressures $ p_i + p_r $. To isolate the pure incident signal, a separate simulation is conducted with the same setup but with a single phase. The reflected signal $ p_r $ is then obtained by subtracting the incident pressure $ p_i $—obtained from the single-phase simulation—from the total pressure measured by the left probe in the two-phase setup. This process enables an accurate inference of the reflected pressure wave, as will be demonstrated in the following sections. 

We can then define the transmission coefficient $ T $ and the reflection coefficient $ R $ as the ratios of the amplitudes of the transmitted wave $ p_t $ and the reflected wave $ p_r $ to the incident wave $ p_i $. Specifically, these are defined as:
\begin{equation}
    T_{\text{num}} := \frac{\max(p_t)}{\max(p_i)}, \qquad R_{\text{num}} := \frac{\max(p_r)}{\max(p_i)},
\end{equation}

and are referred to as the numerical transmission and reflection coefficients, respectively. The subscript "num" indicates that these values are numerical and are calculated from the simulations. 

Table~\ref{tab:1Dparameters} reports the parameters used in all subsequent simulations unless otherwise mentioned. 
The choice of these parameters is primarily driven by the acoustic transmission characteristics between air and water at audible frequencies. This motivation guides the selection of sound speeds in both media and restricts the frequency range considered to values below or on the order of tens of $ \, \SI{}{\kilo\hertz}$. Additionally, we aim to exclude any potential effects caused by surface tension. Regardless, for frequencies in this range, surface tension and viscosity have little to no influence on acoustic transmission and reflection signals when considering air to water transmission \cite{krechetnikov2019viscosity}. 

\begin{table}[h]
\centering
\resizebox{\textwidth}{!}{
\begin{tabular}{@{}ccccccccccc@{}}
\toprule
\multicolumn{3}{c}{Fluid 1} & \multicolumn{3}{c}{Fluid 2} & \multicolumn{3}{c}{Interface} & \multicolumn{1}{c}{Source} \\ \midrule
 $\rho_1 \, (\SI{}{\kilogram\per\meter^3})$  & $c_{s1} \, (\SI{}{\meter\per\second})$  & \multicolumn{1}{c|}{$\eta_1 \, (\SI{}{\pascal\second})$}   & $\rho_2 $  & $c_{s2} $  & \multicolumn{1}{c|}{$\eta_2$}    & $\varepsilon \, (\SI{}{\meter})$ & $\sigma \, (\SI{}{\newton\per\meter})$   & \multicolumn{1}{c|}{$M_0 \, (\SI{}{\meter^3 \second \per \kilogram})$} & $f \, (\SI{}{\hertz})$  \\ \midrule
 $1$ & $343$ & \multicolumn{1}{c|}{$10^{-16}$} & $2$  & $1481$ & \multicolumn{1}{c|}{$10^{-16}$} & $0.01$ & $10^{-16}$ & \multicolumn{1}{c|}{$0.01$} & $1000$       \\ \bottomrule
\end{tabular}
}
\caption{Simulation parameters.}
\label{tab:1Dparameters}
\end{table}

\subsubsection{Spectral convergence}
This section demonstrates the spectral convergence of our implementation by evaluating numerical errors in the transmission and reflection coefficients across different mesh resolutions with varying polynomial orders. To this end, a series of numerical simulations are conducted on a one-dimensional mesh with an increasing number of uniform elements, starting from $N_{el} = 125$, which are then subsequently doubled until $N_{el} = 2000$. The polynomial order is varied from $N=1$ to $N=6$. 

Absolute errors are calculated as the difference between the reference and computed numerical solutions:
\begin{equation}
     \lvert R_{\text{num,best}} - R_{\text{num}} \rvert, \qquad 
     \lvert T_{\text{num,best}} - T_{\text{num}} \rvert, 
     \label{eq::e_num}
\end{equation}
for the errors in the reflection and transmission coefficients. The reference values are obtained from the results of the most accurate numerical solution provided by the finest mesh and a polynomial order $N=7$. The corresponding reference transmission and reflection coefficients are reported in Table~\ref{tab:reference_coefficients}.
\begin{table}[h]
\centering
\begin{tabular}{@{}ll@{}}
\toprule
Reference value for the transmission coefficient $ R_{\text{num,best}} $                  & $0.790988459995992$  \\ \midrule
Reference value for the reflection coefficient $ T_{\text{num,best}} $ & $1.797929864428584$ \\ \bottomrule
\end{tabular}
\caption{Transmission and reflection coefficients calculated from the mesh $N_{el}=2000$ and polynomial order $N=7$.}
\label{tab:reference_coefficients}
\end{table}
 
\begin{figure}[ht]
    \centering
    \begin{subfigure}{0.48\textwidth}
        \centering
        \includegraphics[width=\linewidth]{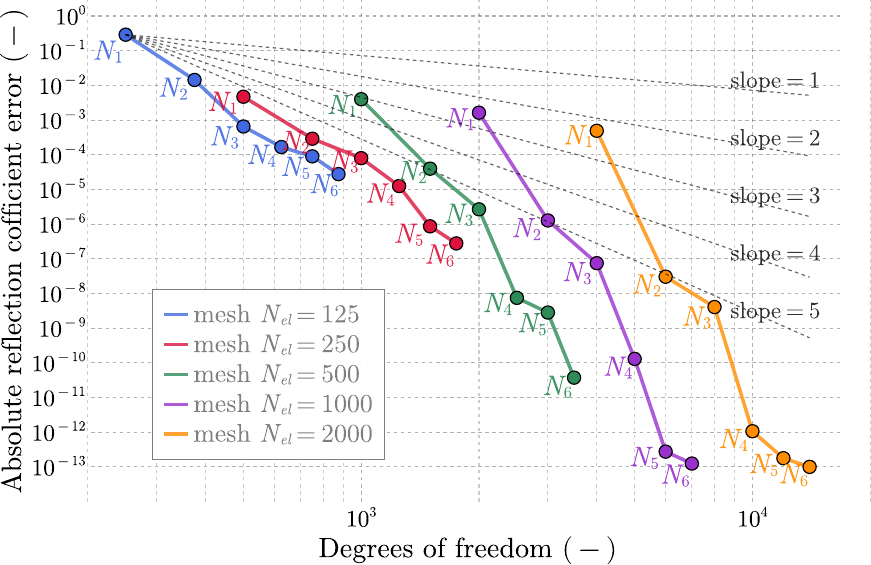}
        \caption{}
        \label{fig:1D1_R_errors}
    \end{subfigure}
    \hspace{0.002\textwidth} 
    \begin{subfigure}{0.48\textwidth}
        \centering
        \includegraphics[width=\linewidth]{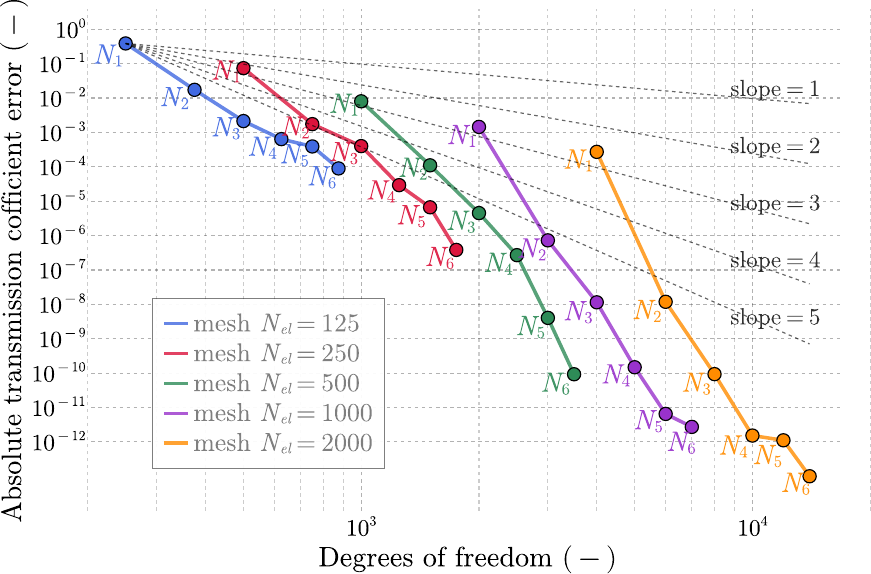}
        \caption{}
        \label{fig:1D1_T_errors}
    \end{subfigure}
    \caption{Absolute errors for reflection (\subref{fig:1D1_R_errors}), and transmission (\subref{fig:1D1_T_errors}), coefficients vs degrees of freedom given by different meshes and polynomial orders. Errors decrease with increased resolution}
    \label{fig:1d1_errors}
\end{figure}

Fig.~\ref{fig:1d1_errors} presents errors reduced to almost machine precision in both reflection (Fig.~\ref{fig:1D1_R_errors}) and transmission (Fig.~\ref{fig:1D1_T_errors}) as a function of the number of degrees of freedom (DOFs), calculated as: $\text{DOFs} = N_{el} \times (N+1)$. As expected from a high-order numerical method, the solver exhibits spectral convergence, with errors decreasing exponentially as the polynomial order increases. Furthermore, for a fixed number of DOFs, a relatively coarse mesh with a higher polynomial order achieves greater accuracy compared to a finer mesh with a lower polynomial order, highlighting the advantages of high-order discretization.

\subsubsection{Modeling errors}
Although the inclusion of a diffuse interface enhances robustness and simplifies the implementation of a numerical scheme by smearing the interface across several elements, as opposed to using a sharp interface, it is expected that a diffuse interface, such as that provided by the Cahn-Hilliard model, introduces modeling errors. In the previous section, we analyzed the numerical errors and confirmed the expected spectral convergence. In this section, we focus on quantifying the modeling errors for acoustic transmission and reflection introduced by approximating the sharp interface with a diffuse interface. To achieve this, we evaluate the modeling errors as $\varepsilon \to 0$, considering analytical expressions for wave transmission and reflection derived for a sharp interface. We use the same setup as in the previous section for our simulations, with identical parameters for both fluids.

For linear plane waves with normal incidence, the analytical expressions for the reflection and transmission coefficients are given by \cite{brekhovskikh2012acoustics}:
\begin{equation}
    R_\text{exact} = \frac{z_2 - z_1}{z_1+z_2}, \qquad
    T_\text{exact} = \frac{2 z_2 }{z_1+z_2},  
\label{eq:planeTR}
\end{equation}
where $z_1$ and $z_2$ are the impedance coefficients of each phase and are given by the product $z_i=\rho_i c_{si},  i \in \{0,1\}$. We aim to monitor the convergence of the numerical transmission and reflection coefficients toward those specified in equation \eqref{eq:planeTR} as we vary the thickness of the interface, varying $\varepsilon$. To isolate the effects of modeling errors, we minimize numerical errors by selecting the mesh with $N_{el} = 500$ and a polynomial order $N=6$ for all simulations. As demonstrated in the previous section, this level of spatial accuracy achieves sufficiently small discretization errors, effectively nullifying their contribution. The absolute modeling errors are calculated as:
\begin{equation}
     \lvert R_{\text{exact}} - R_{\text{num}} \rvert, \qquad
     \lvert T_{\text{exact}} - T_{\text{num}} \rvert,
     \label{eq::e_mod}
\end{equation}
 for modeling errors in the reflection and transmission coefficients.  
\begin{figure}[h] 
    \centering
    \includegraphics[width=0.6\textwidth]{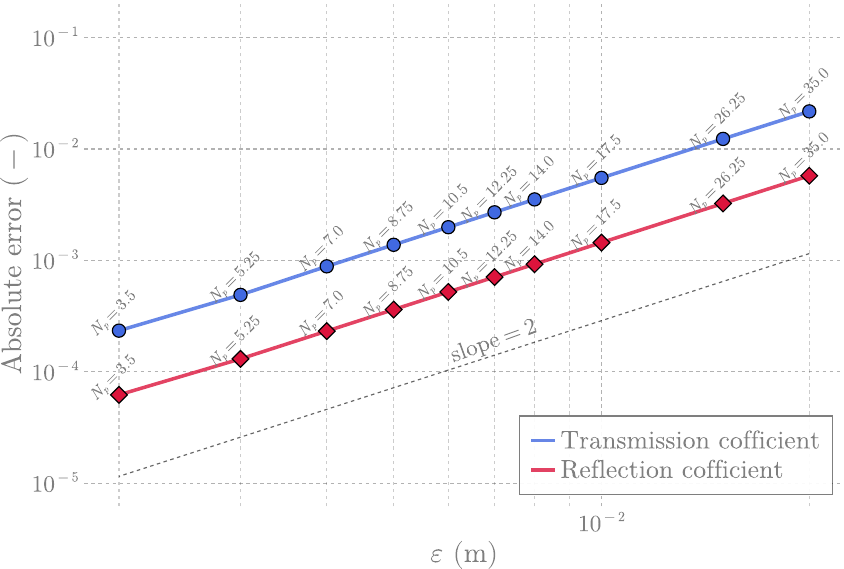}
    \caption{Transmission and reflection modeling errors versus $\epsilon$ at $f = \SI{1000}{\hertz}$. Errors decay with $\varepsilon^2$.}
    \label{fig:1D3_errorvsepsi}
\end{figure}

Fig.~\ref{fig:1D3_errorvsepsi} shows the modeling errors in the transmission and reflection coefficients as a function of half the interface thickness, $\varepsilon \in [0.002, 0.02]$, for a set frequency of $f = \SI{1000}{\hertz}$. The number of DOFs within the interface, $N_p$, is reported for each point in Fig.~\ref{fig:1D3_errorvsepsi}, ensuring that there are enough points within the interface. $N_p$ is approximated to be:
\begin{equation}
    N_p = 2  \varepsilon   N_{el} (N+1) / L ,
    \label{eq:Npoints}
\end{equation}
 where $L$ is the width of the domain. It can be seen that the modeling errors are second order with respect to the interface thickness, i.e.,  $\sim O(\varepsilon^2)$.

\begin{figure}[ht]
    \centering
    \begin{subfigure}{0.47\textwidth}
        \centering
        \includegraphics[width=\linewidth]{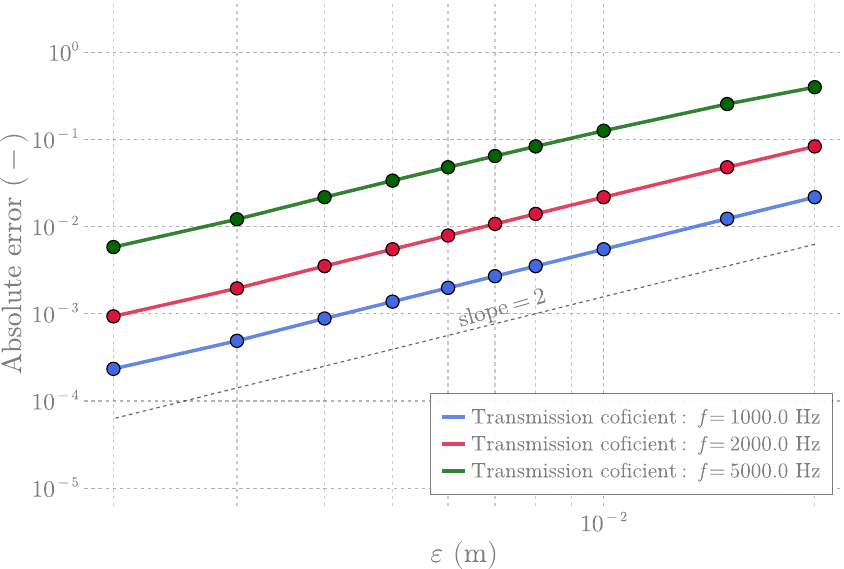}
        \caption{}
        \label{fig:1D4_errorvsepsiandF}
    \end{subfigure}
    \hspace{0.002\textwidth} 
    \begin{subfigure}{0.47\textwidth}
        \centering
        \includegraphics[width=\linewidth]{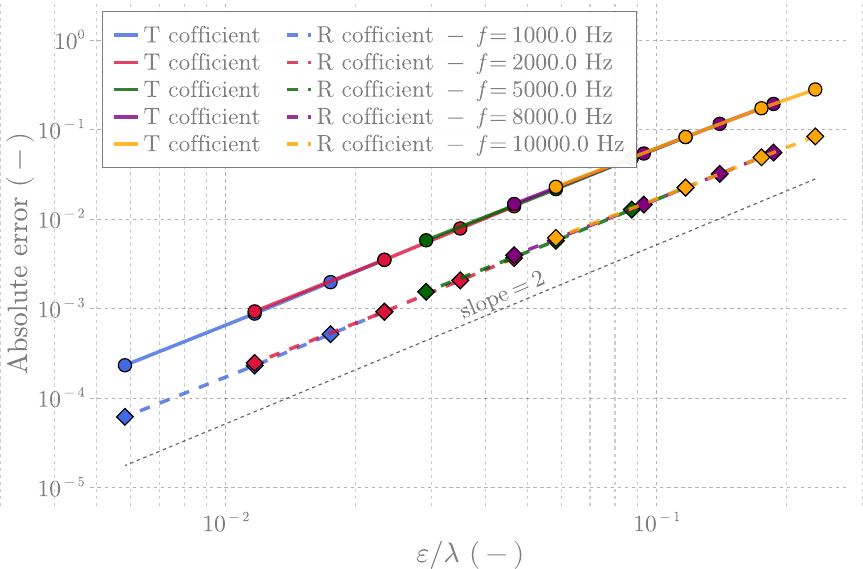}
        \caption{}
        \label{fig:1D5_errorvsepsiandlambda}
    \end{subfigure}
    \caption{Modeling errors in (a) transmission versus $\varepsilon$ for frequencies 1000, 2000, and 5000 Hz, and (b) transmission and reflection versus $\varepsilon$ scaled by the wavelength $\lambda$ for frequencies 1000, 2000, 5000, 8000, and 10000 Hz. Errors decay as $\varepsilon^2$.}
\end{figure}

To assess the impact of the incident signal frequency on the error, we repeated the previous simulation using different frequencies. Fig.~\ref{fig:1D4_errorvsepsiandF} shows the modeling transmission error as a function of $\varepsilon$ for different frequencies. The reflection error is not shown here for brevity, but also displays similar behavior. It is evident that the errors increase with increasing frequencies but maintain the $ \sim O(\varepsilon^2)$ rate. 

Fig.~\ref{fig:1D5_errorvsepsiandlambda} considers the same data but scales the interface thickness with the wavelength $\lambda =c_{s1} / f$ on the x-axis and considers frequencies up to $\SI{10}{\kilo\hertz}$. As shown in Fig.~\ref{fig:1D5_errorvsepsiandlambda}, all errors are superimposed for both reflection and transmission errors. This implies that modeling errors actually scale with $O\left(\varepsilon / \lambda\right)^2$. The convergence rates we measure for the sharp-interface limit align with previous numerical studies of the incompressible CH-NS system (i.e., without acoustics), which report orders of $ O(\varepsilon^\alpha)$ where $ 0<\alpha<3 $. Our results show that as $\varepsilon \to 0$, the model approaches the sharp-interface limit for acoustic propagation across the diffuse interface.

\subsubsection{Considering both numerical and modeling errors}
In the previous sections, we isolated both numerical and modeling errors. We now present the results that consider both sources of error concurrently for the transmission coefficient. Similar results were obtained for the reflection coefficient, but were omitted for brevity.

\begin{figure}[ht]
    \centering
    \begin{subfigure}{0.47\textwidth}
        \centering
        \includegraphics[width=\linewidth]{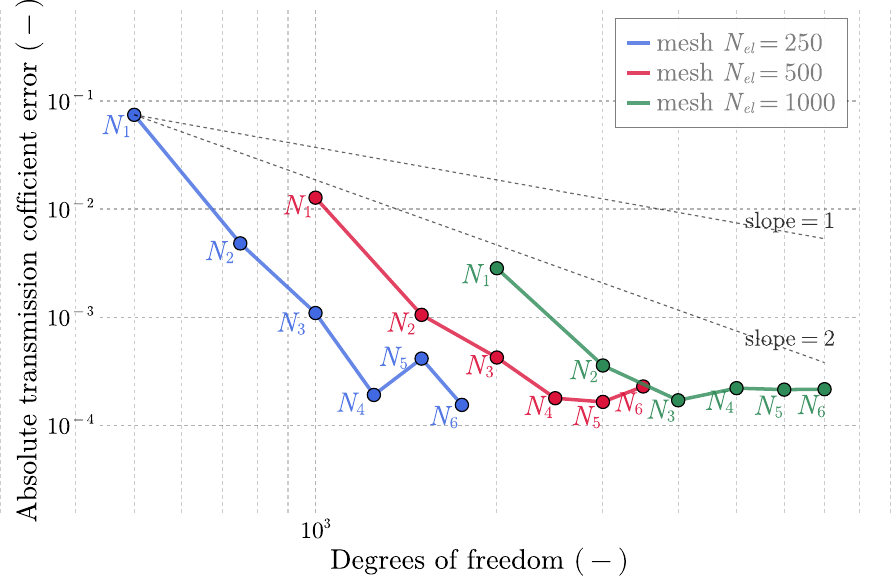}
        \caption{}
        \label{fig:1D6_transmission_error}
    \end{subfigure}
    \hspace{0.002\textwidth} 
    \begin{subfigure}{0.47\textwidth}
        \centering
        \includegraphics[width=\linewidth]{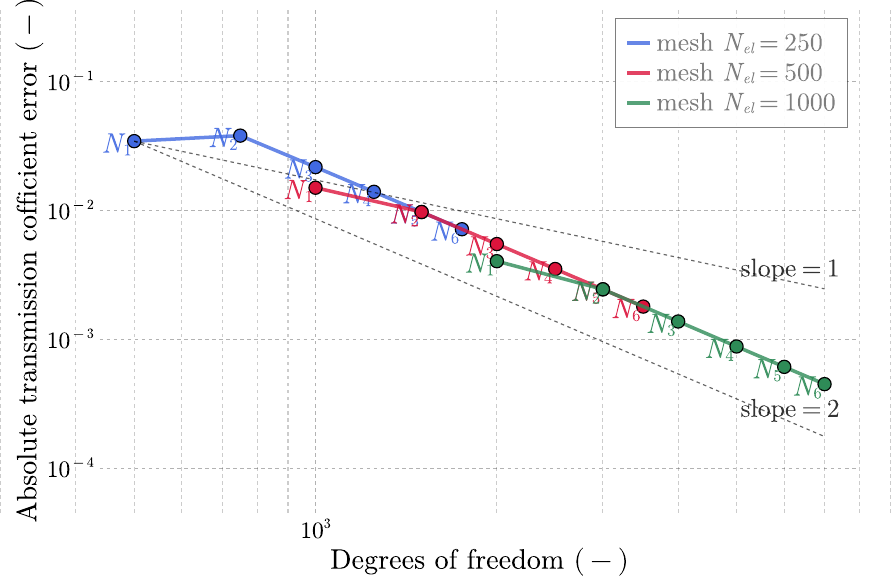}
        \caption{}
        \label{fig:1D8_transmission_error}
    \end{subfigure}
\caption{Absolute transmission errors including numerical and modeling contributions. (\subref{fig:1D6_transmission_error}) Errors for fixed $\varepsilon = 2 \times 10^{-3}$, which decay exponentially before saturating. (\subref{fig:1D8_transmission_error}) Errors with $\varepsilon$ scaled to have $10 N_p$ points within the interface, showing quadratic decay}
    \label{fig:onepsi_errors}
\end{figure}

Following our methodology from the initial numerical experiments, we perform the same simulations but compute the error with respect to the analytical solution, as in the previous section using \eqref{eq::e_mod}. The results are reported in Fig.~\ref{fig:1D6_transmission_error}, where we selected an interface thickness of $\varepsilon = 2 \times 10^{-3}$ and a frequency of $f = \SI{1000}{\hertz}$. The curves show a high order of convergence initially, before reaching a plateau at approximately an error of $2.2 \times 10^{-4}$. Referring to Fig.~\ref{fig:1D4_errorvsepsiandF} ($f = \SI{1000}{\hertz}$, $\varepsilon=2 \times 10^{-3}$), we infer a similar error, which implies that the plateau is due to modeling errors. 

We repeat the same set of simulations, but rather than maintaining a constant interface thickness, we scale $\varepsilon$ so that there are 10 $N_p$ points within the interface according to equation \eqref{eq:Npoints}. Ensuring that at least 10 points are placed within the interface aligns with established practices \cite{boyanova2014efficient,teigen2011diffuse}. Fig.~\ref{fig:1D8_transmission_error} reports the results and illustrates how the convergence rates revert to second order, dominated by modeling errors.

In these experiments, we have established that the total error $E_\text{tot}$ is the sum of numerical $E_\text{num}$, and modeling $E_\text{model}$ errors:
\begin{equation}
    E = E_{num} + E_{model}.
\end{equation}
Taking into account both error contributions, while $E_\text{num}$ decays exponentially with increasing polynomial order, $E_\text{model}$ decays quadratically, dominating the contribution to total error, which explains the spectral convergence followed by the second-order rate in Fig.~\ref{fig:1D8_transmission_error}. However, values ultimately converge to the sharp-interface limit for analytical acoustic expressions for transmission and reflection. 

\subsection{Two dimensions - Snell's law}
\label{subsec:2D}
When an acoustic wave encounters an interface separating two media with distinct sound speeds at an oblique angle, the transmitted wave undergoes a change in direction, a phenomenon known as wave refraction. This phenomenon is governed by Snell's law, which relates the angles of incidence $\theta_i$ and transmission $\theta_t$ to the respective sound speeds in each medium:
\begin{equation}
\frac{\sin(\theta_i)}{c_{s1}} = \frac{\sin(\theta_t)}{c_{s2}}.
\label{eq:Snell}
\end{equation}  

To assess the accuracy of our numerical model in capturing this fundamental wave behavior, we extend our prior one-dimensional analysis to a two-dimensional setup. The incident plane wave remains excited by the source term of equation \eqref{eq:source}, but the interface is slanted at an angle, which in this case is the same angle of incidence $\theta_i$. Thus, the initial concentration profile is given by:
\begin{equation}
    c_0 = 1.0 - 0.5\left(1+tanh \left( 2\dfrac {cos(\theta_i)x+sin(\theta_i)y } {\varepsilon} \right)\right).
    \label{eq:IC_Angle}
\end{equation}
The mesh is 2D cartesian, so the total number of degrees of freedom now is $\text{DOFs} = N_{el}\times (N+1)^2$, for a polynomial order $N$.
The upper and lower boundary conditions are designated as impermeable slip and no flux for concentration and chemical potential, as stated by \eqref{eq:no-flux}.

\begin{figure}[ht]
    \centering
    \begin{subfigure}{0.47\textwidth}
        \centering
        \includegraphics[width=\linewidth]{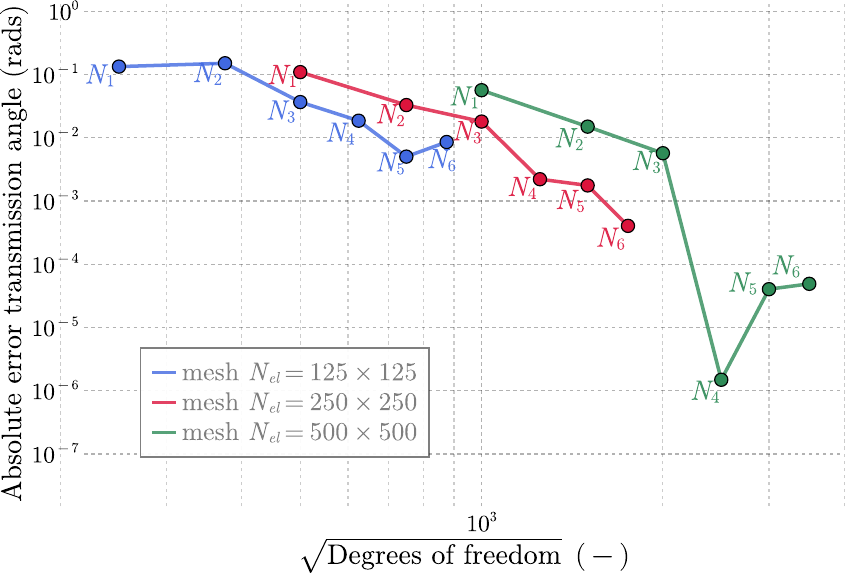}
        \caption{}
        \label{fig:2D2_vsDOFs10}
    \end{subfigure}
    \hspace{0.002\textwidth} 
    \begin{subfigure}{0.47\textwidth}
        \centering
        \includegraphics[width=\linewidth]{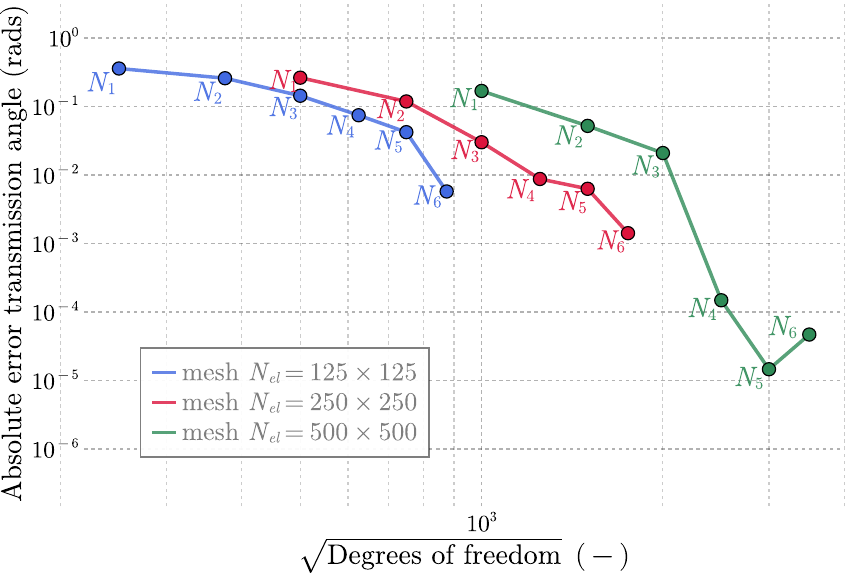}
        \caption{}
        \label{fig:2D2_vsDOFs13}
    \end{subfigure}
\caption{Convergence rates for the absolute error in the transmission angle. (\subref{fig:2D2_vsDOFs10}) $\theta_i = 10^\circ$ and (\subref{fig:2D2_vsDOFs13}) $\theta_i = 13^\circ$, plotted as a function of the square root of the degrees of freedom, showing spectral convergence.}
    \label{fig:2D2vsDOFs}
\end{figure}

A series of numerical experiments are performed to validate the solver's ability to accurately reproduce Snell's law. The chosen parameters align with those used in the one-dimensional test cases outlined in Table~\ref{tab:1Dparameters}, except that now we consider a higher wave frequency of $\SI{10}{\kilo\hertz}$. For the selected sound speed values, the critical incidence angle is computed as:
\begin{equation*}
\theta_c = \arcsin\left(\dfrac{c_{s1}}{c_{s2}} \right) = \arcsin\left(\dfrac{343}{1481} \right) \approx 13.39^\circ,
\end{equation*}
beyond which, incident waves are meant to experience total reflection.

Fig.~\ref{fig:2D2vsDOFs} shows the spectral convergence of the solver for an incident angle of $\theta_i=10^\circ$ (Fig.~\ref{fig:2D2_vsDOFs10}), and $\theta_i=13^\circ$ (Fig.~\ref{fig:2D2_vsDOFs13}) where the latter angle is close to the critical angle. The results confirm that the numerical solution aligns with Snell's law as the number of degrees of freedom increases. Errors are computed relative to the analytical expression in equation \eqref{eq:Snell}, and are reported in radians. More details on how the transmission angle was computed are included in \ref{sec:Appendix-A}.

\begin{figure}[ht]
    \centering
    \begin{subfigure}{0.47\textwidth}
        \centering
        \includegraphics[width=\linewidth]{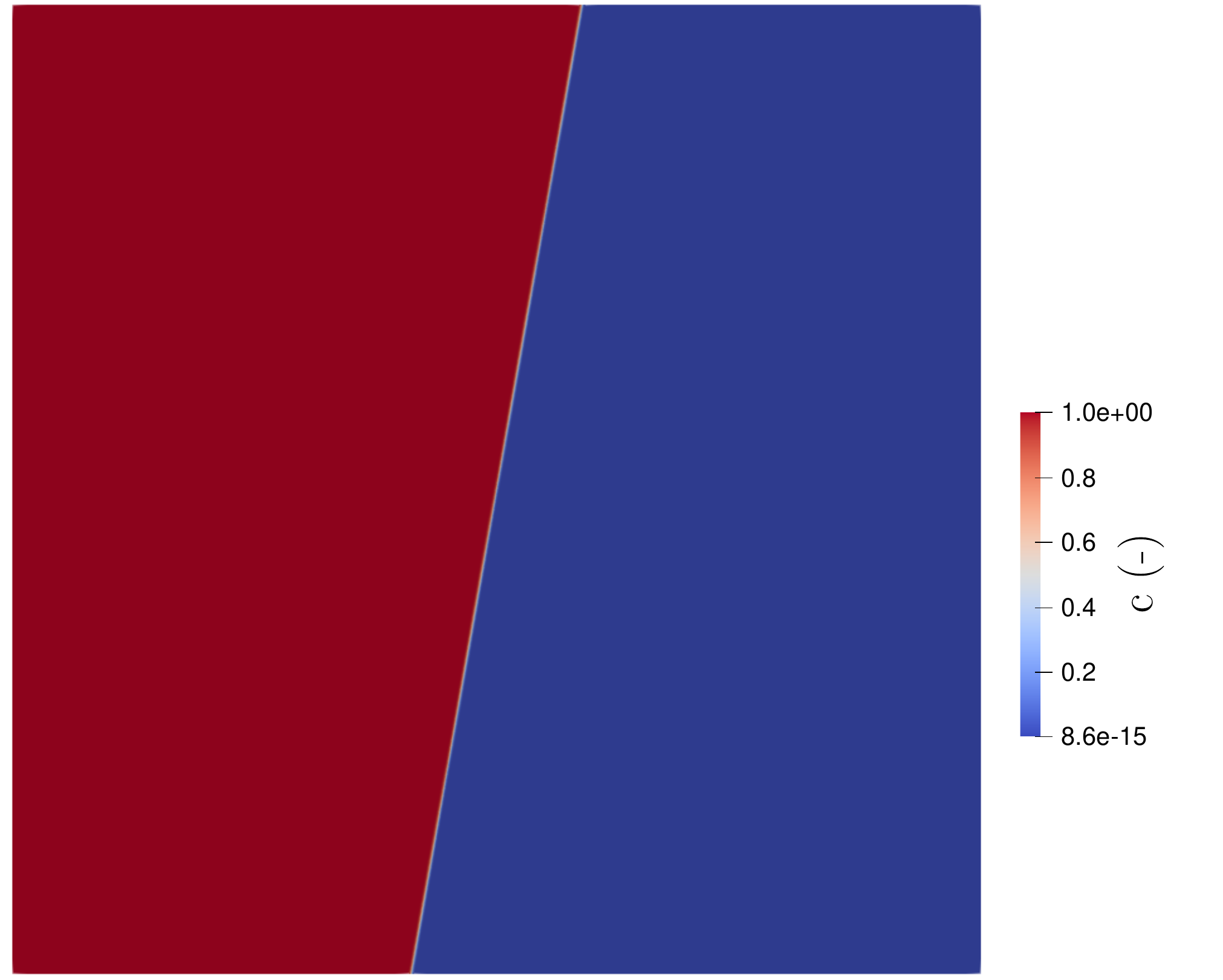}
        \caption{}
        \label{fig:Snell_c}
    \end{subfigure}
    \hspace{0.002\textwidth} 
    \begin{subfigure}{0.47\textwidth}
        \centering
        \includegraphics[width=\linewidth]{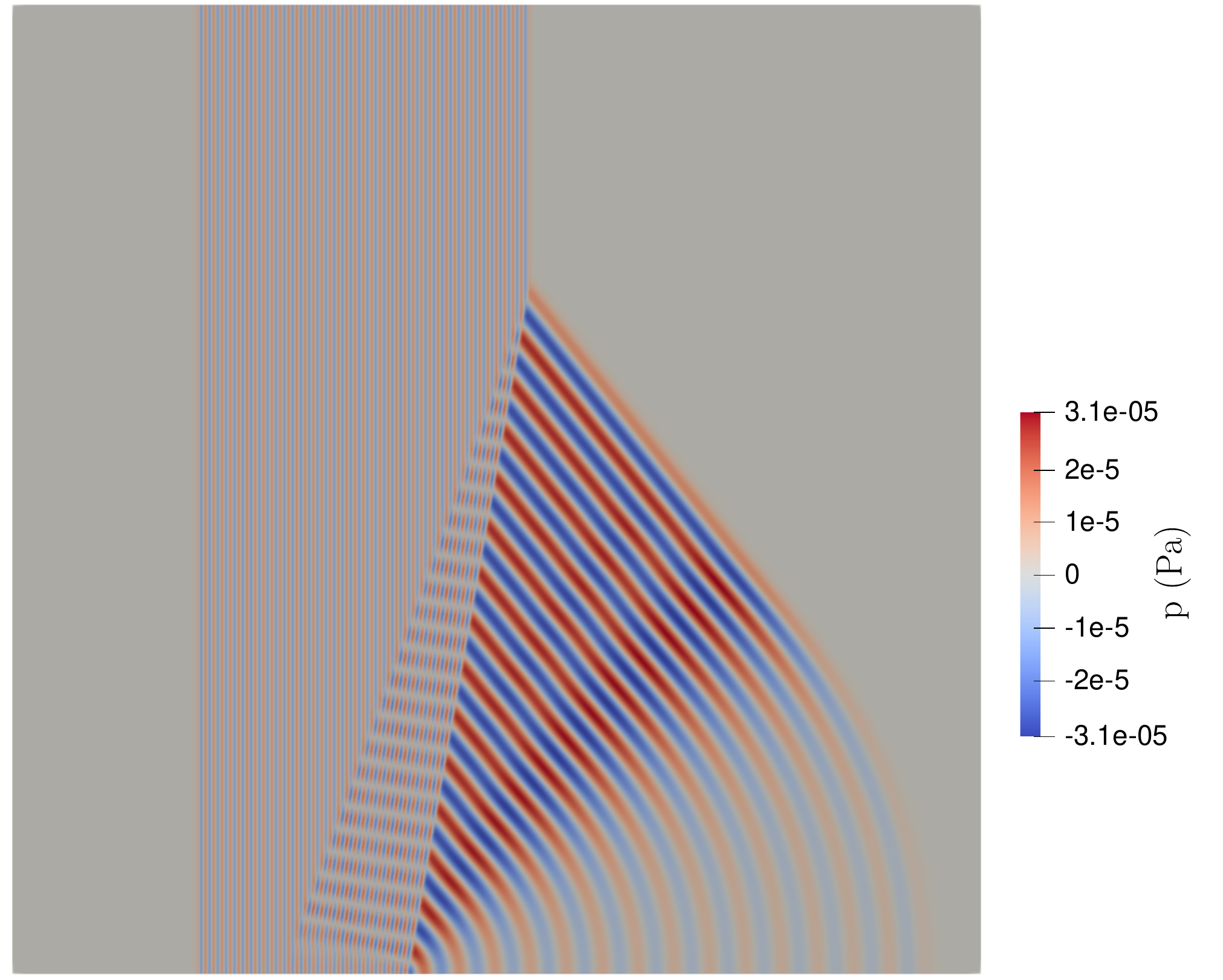}
        \caption{}
        \label{fig:Snell_p}
    \end{subfigure}
\caption{Concentration (\subref{fig:Snell_c}) and pressure (\subref{fig:Snell_p}) fields at the end of the simulation for $\theta_i = 10^\circ$ and $f = \SI{10}{\kilo\hertz}$, showing refraction of the wave at the interface.}
    \label{fig:SNELL}
\end{figure}

Fig.~\ref{fig:Snell_c} shows the concentration field $\theta_i=10^\circ$, while Fig.~\ref{fig:Snell_p} shows the pressure at the end of the simulation for the best setup. The pressure contour plots exhibit good qualitative agreement with the expected behavior of wave refraction.  

\begin{figure}[ht]
\centering
\includegraphics[width=0.6\linewidth]{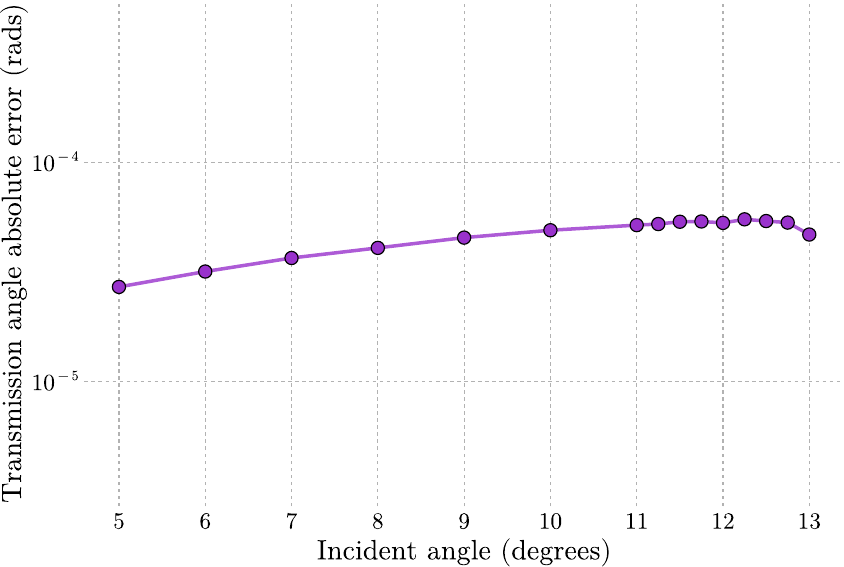}
\caption{Error variation of the transmitted angle with varying angle of incidence, showing minimal error for all incident angles.}
\label{fig:2D4_vsalpha}
\end{figure}

To examine the effect of the angle of incidence on numerical errors, additional simulations are performed by varying $\theta_i$. The finest 2D setup from the previous experiment, with $N_{el}=500 \times 500$ and polynomial order $N=6$, yielding over 12 million DOFs, is selected. Fig.~\ref{fig:2D4_vsalpha} illustrates that errors increase as the incident wave approaches the critical angle. Nevertheless, errors remain small across all incident angles, demonstrating the solver's capability to accurately reproduce Snell's law.

\subsection{Three dimensions - Spherical wave}
\label{subsec:3D}
In this final numerical experiment, we investigate the propagation of three-dimensional spherical waves through two media separated by a flat surface. The setup here is similar to that in the previous sections. The computational domain is defined by a cube with a side length of $2 \SI{}{\meter}$, the fluid interface being placed at the $z=0$ plane. A point source is positioned at $(0,\,0.5,\,0)$ above the interface. The structured mesh includes a locally refined region that encompasses four elements across the interface, with progressive grading to increase resolution near the interface. The simulation employs polynomial order $N = 4$ and a number of elements $N_e=540,000$, yielding a total of 
$67,500,000$ DOFs. A sufficient number of $20$ DOFs were placed along the interface. The physical parameters remain the same as in the previous sections, but here we chose a high density ratio by setting $\rho_2=1000 \SI{}{\kilogram\per\meter^3}$ to mimic air-water conditions. 

Regarding the selection of $\varepsilon$, we performed an analysis analogous to the one conducted in the one-dimensional case, but set $\rho_2=1000 \SI{}{\kilogram\per\meter^3}$. This allowed us to obtain trends similar to those shown in Fig.~\ref{fig:1D5_errorvsepsiandlambda}. From these results, the following empirical modeling errors for the transmission $E_{T,\text{model}}$ and reflection $E_{R,\text{model}}$ coefficients can be derived:
\begin{equation}
    E_{T,\text{model}} = 54.88 \left(\frac{\varepsilon}{\lambda}\right)^2, \quad 
    E_{R,\text{model}} = 0.3077 \left(\frac{\varepsilon}{\lambda}\right)^2.
\end{equation}
Assuming an accuracy requirement of $10^{-3}$, the value of $\varepsilon$ can be estimated from $E_{T,\text{model}}$, as it imposes a more stringent constraint compared to $E_{R,\text{model}}$, leading to $\varepsilon \approx 0.00145~\text{m}$

\begin{figure}[ht]
\centering
\includegraphics[width=0.7\linewidth]{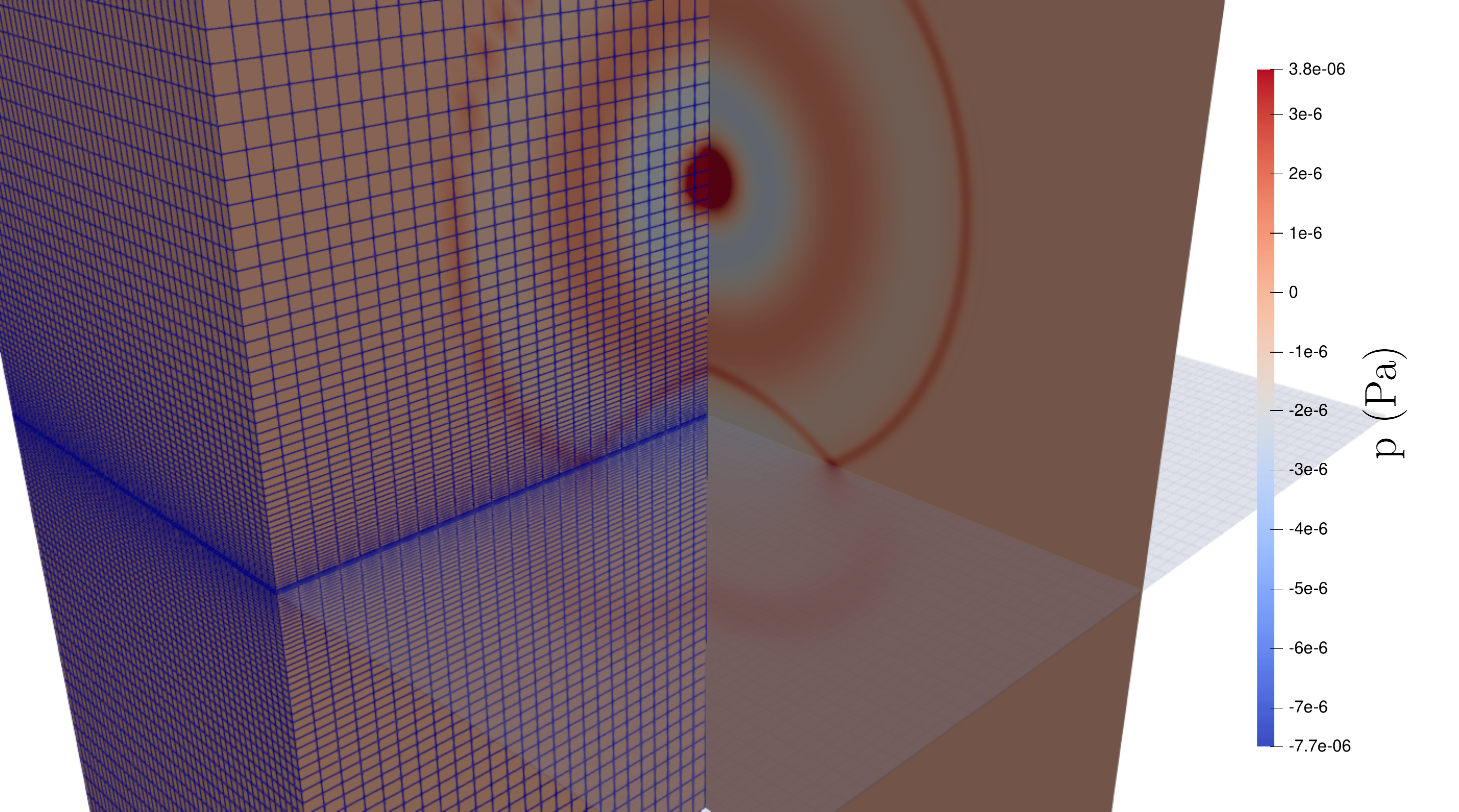}
\caption{Pressure contours on the $xy$ and $yz$ planes with the computational mesh at $t = 1.75 \times 10^{-3} \SI{}{\second}$. The contours illustrate the wave propagation across the two media, showing aa transmitted and a reflected wave at the phase interface.}
\label{fig:3D}
\end{figure}

Fig.~\ref{fig:3D} shows two orthogonal plane slices of the domain, depicting the geometry, mesh, and the pressure waves at $t = 1.75 \times 10^{-3} \SI{}{\second}$ right after passing into the second phase. A reflected wave is also visible. Mesh refinement is applied at the interface between the two phases, with a transparent yz-plane. shown. 

\begin{figure}[ht]
\centering
\includegraphics[width=0.7\linewidth]{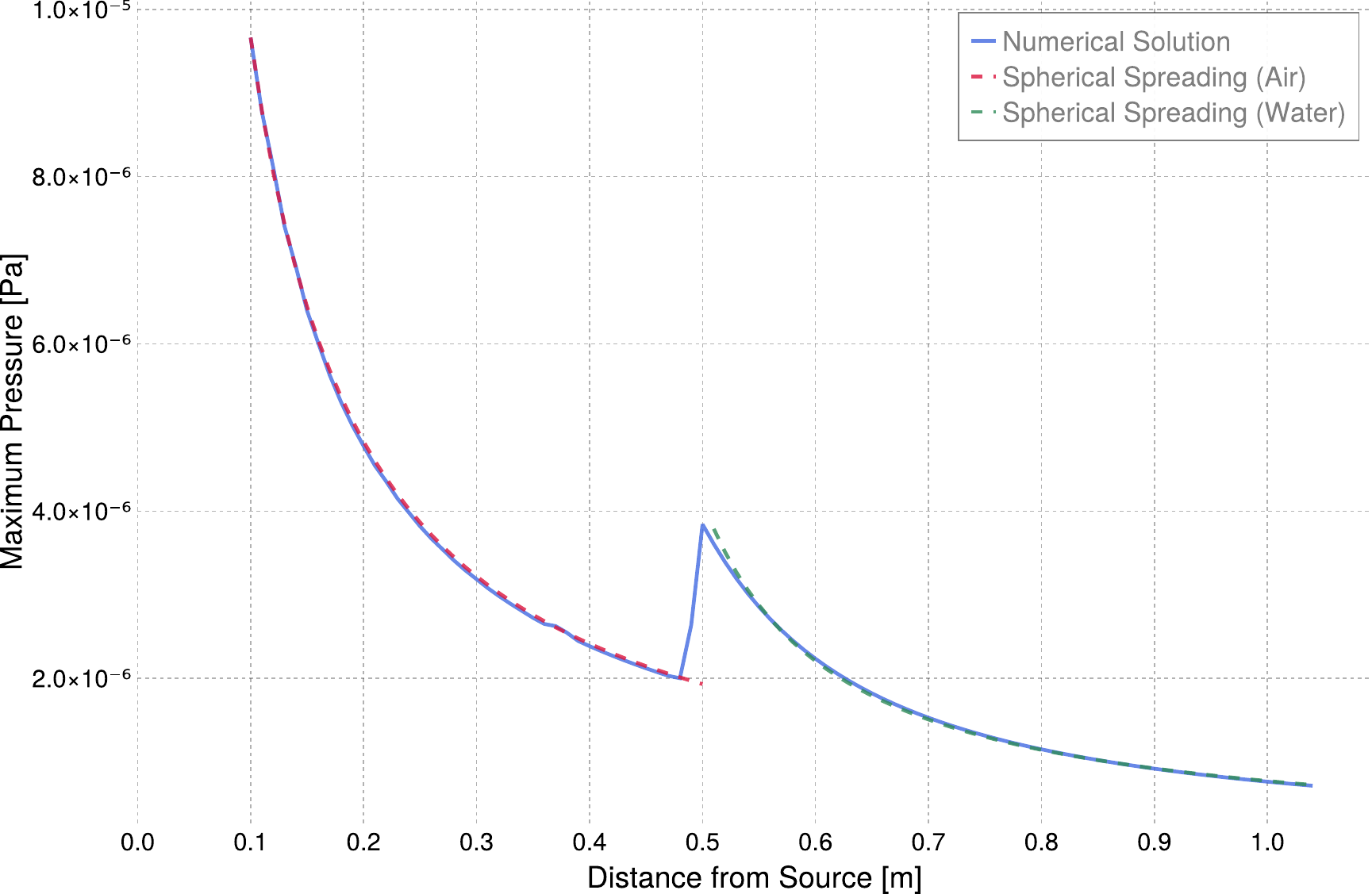}
\caption{Pressure wave transmission from air to water following the theoretical spherical spreading.}
\label{fig:3Dspreading}
\end{figure}

Fig.~\ref{fig:3Dspreading} shows the maximum pressure recorded at a given distance from the source. The graph exhibits three distinct stages for wave propagation. In the first stage, the pressure amplitude decays as $1/r$, consistent with the theoretical behavior of a spherical wave, where $r$ denotes the distance from the source. As the wave crosses into the second phase, a sharp increase in pressure occurs, in line with the prediction of the transmission equation~\eqref{eq:Cahn-Hilliard}. Finally, the pressure wave decays proportionally to $1/r_{\text{eff}}$, as it propagates through the second phase, where $r_{\text{eff}}$ is defined as the distance to a source placed at $h c_{s1}/c_{s2}$ with $h$ being the distance from the source to the interface ~\cite[Sec.~8.7]{pierce2019acoustics}. Fig.~\ref{fig:3Dspreading} shows excellent agreement between the numerical and theoretical results in all three stages. 

\section{Conclusion}\label{sec:Conclusion}
This work presents a novel approach to modeling acoustic wave propagation in multiphase media using a diffuse interface formulation with weak compressibility. We validated the method through comprehensive numerical experiments in one-, two-, and three-dimensional settings, demonstrating the spectral convergence of the proposed method, comparing it with analytical solutions for acoustic transmission and reflection coefficients, quantifying the modeling errors introduced by the diffuse interface width parameter, qualitatively and quantitatively validating Snell's law, and showcasing the ability to simulate high density ratios. Our analysis quantifies both numerical and modeling errors, confirming that the method asymptotically approaches the sharp-interface limit. The results highlight the potential of diffuse interface models for high-fidelity acoustic simulations in multiphase environments. This can prove to be useful for direct acoustic computation in incompressible multiphase flows. The approach preserves the benefits of diffuse interface models, such as the natural incorporation of surface tension effects and the ability to naturally reproduce and track the interface without the need for any explicit interface handling, while accurately capturing acoustic transmission and reflection phenomena. Future work will focus on exploring marine aero/hydro-acoustics.

\section*{Acknowledgments}
AB, GN, GR and EF acknowledge the funding received by the Grant DeepCFD (Project No. PID2022-137899OB-I00) funded by MICIU/AEI/10.13039/501100011033 and by ERDF, EU. 
OM, and EF acknowledge the funding from the European Union (ERC, Off-coustics, project number 101086075). Views and opinions expressed are, however, those of the authors only and do not necessarily reflect those of the European Union or the European Research Council. Neither the European Union nor the granting authority can be held responsible for them.

All authors gratefully acknowledge the Universidad Politécnica de Madrid (www.upm.es) for providing computing resources on Magerit Supercomputer and the computer resources at MareNostrum and the technical support provided by Barcelona Supercomputing Center (projects RES-IM-2024-1-0003 and RES-IM-2025-1-0011).
Finally, the authors gratefully acknowledge the EuroHPC JU for the project EHPC-REG-2023R03-068 for providing computing resources of the HPC system Vega at the Institute of Information Science.

\begin{appendices}
\section{Calculation of transmission}
\label{sec:Appendix-A}

To compute the transmission angle from the two-dimensional simulations, a probe is placed to the right of the interface at the mid-vertical location of the domain and at $x=\SI{0.02}{\meter}$, well into the second phase and away from the interface. The probe monitors the pressure of the passing wave and its velocity components. The numerical angle of transmission can then be inferred directly from the relation:
\begin{equation}
   \theta_{t,num} = \arctan \left( \dfrac{v}{u} \right).
\end{equation}
At the start of the simulation, $u$ is at machine precision zero, as there are no passing waves yet. Therefore, we set a threshold for the pressure value at the probe; below this threshold, the angle of transmission is not computed and is set to zero. The method for calculating the angle now becomes:
\begin{equation}
\theta_{t,num} =
\begin{cases} 
0, & \text{if } p < t_f \cdot p_{\max} \\
\arctan \left( \dfrac{v}{u} \right), & \text{otherwise}.
\end{cases}
\end{equation}
where $p_{max}$ is the maximum value recorded by the probe for the entire time series, and $t_f$ is a threshold factor that we have set to $t_f= 1 \times 10^{-2}$.  

\begin{figure}[ht]
    \centering
    \begin{subfigure}{0.47\textwidth}
        \centering
        \includegraphics[width=\linewidth]{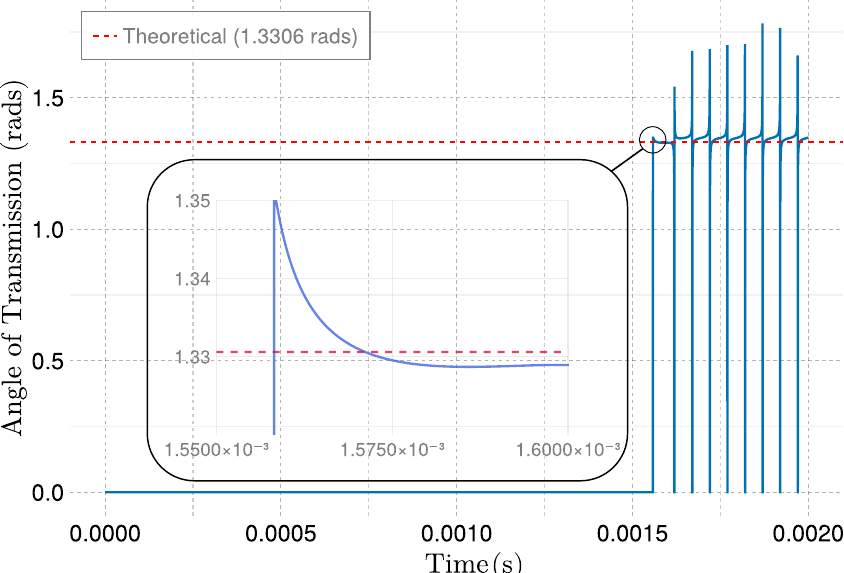}
        \caption{}
        \label{fig:2D2_appendixp3}
    \end{subfigure}
    \hspace{0.002\textwidth} 
    \begin{subfigure}{0.47\textwidth}
        \centering
        \includegraphics[width=\linewidth]{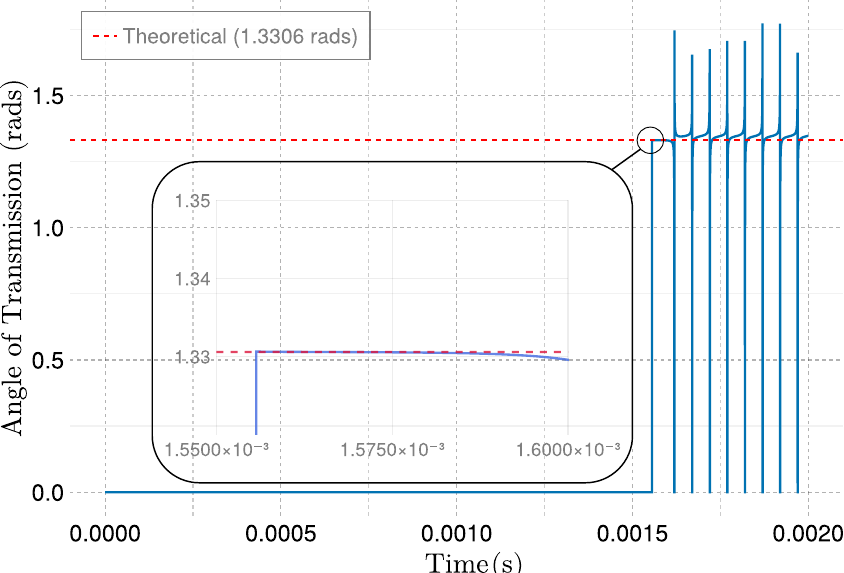}
        \caption{}
        \label{fig:2D2_appendixp6}
    \end{subfigure}
    \caption{Transmission angle $\theta_t$ versus time at the probe located at $x=0.05$ for incidence angle $\theta_i = 13^\circ$, mesh $N_{el} = 500 \times 500$, and frequency $f = \SI{10}{\kilo\hertz}$. (\subref{fig:2D2_appendixp3}) Polynomial order $N=3$, (\subref{fig:2D2_appendixp6}) Polynomial order $N=6$.}
    \label{fig:Annex}
\end{figure}

As an example, Fig.~\ref{fig:Annex} reports the transmission angle calculated by the probe during a simulation for different polynomial orders $N=3$ (Fig.~\ref{fig:2D2_appendixp3}) and $N=6$ (Fig.~\ref{fig:2D2_appendixp6}). The probe initially registers no signal until the transmitted wave passes through it, followed by oscillations due to the signal's oscillatory nature and reflections from the bottom walls. The reported angle of transmission used in the previous sections is the very first one registered by the probe, as shown in the zoomed-in figure.

\section{Discrete entropy balance}
\label{sec:Appendix-B}

The discrete entropy balance can be evaluated following the semi-discrete law derived in \cite{manzanero2020entropyNSCH} (Eqs.~149--150):
\begin{equation}
\frac{d}{dt}
\Bigg(
\bar{\mathcal{E}}^\beta
+
\sum_{f\in\partial\Omega}
\int\limits_{\partial E ,N} f_w(c)\,dS_{\xi}
\Bigg)
=
-\sum_{e \in \Omega}
\big\langle
J\,\big(M_0 |\nabla_\xi \mu|^2 + 2\eta \boldsymbol{S} : \boldsymbol{S}\big), 1
\big\rangle_{E,N}
\le 0,
\label{eq:entropy_balance}
\end{equation}
where \(f_w(c)\) is the boundary free energy function, and the discrete entropy is:
\begin{equation}
\bar{\mathcal{E}}^\beta =
\sum_{e \in \Omega} \langle J\, \mathcal{E} ,1\rangle_{E,N}
+ \frac{3}{4}\sigma \varepsilon
\sum_{f \in \text{int}} \int_{\partial E,N} \beta \, (\jump{c})^2 \, dS_\xi,
\end{equation} 
with element-wise entropy:
\begin{equation}
\mathcal{E} = f_0(c)
+ \frac{3}{4}\sigma \varepsilon |\nabla c|^2
+ \frac{1}{2}\rho |\boldsymbol{u}|^2
+ \frac{p^2}{2 \rho c_s^2}.
\end{equation}

Now, the entropy residual is defined as:
\begin{equation}
\mathcal{R}(\mathcal{E}) = 
\frac{d}{dt}
\Bigg(
\bar{\mathcal{E}}^\beta
+
\sum_{f\in\partial\Omega}
\int_{\partial E,N} f_w(c)\,dS_\xi
\Bigg)
- \sum_{e \in \Omega}
\big\langle
J\,\big(M_0 |\nabla_\xi \mu|^2 + 2 \eta \boldsymbol{S} : \boldsymbol{S} \big), 1
\big\rangle_{E,N},
\label{eq:entropy_residual}
\end{equation}
which must satisfy \(\mathcal{R}(\mathcal{E}) \le 0\) for an entropy-stable scheme. 

\begin{figure}[ht]
\centering
\includegraphics[width=1.0\linewidth]{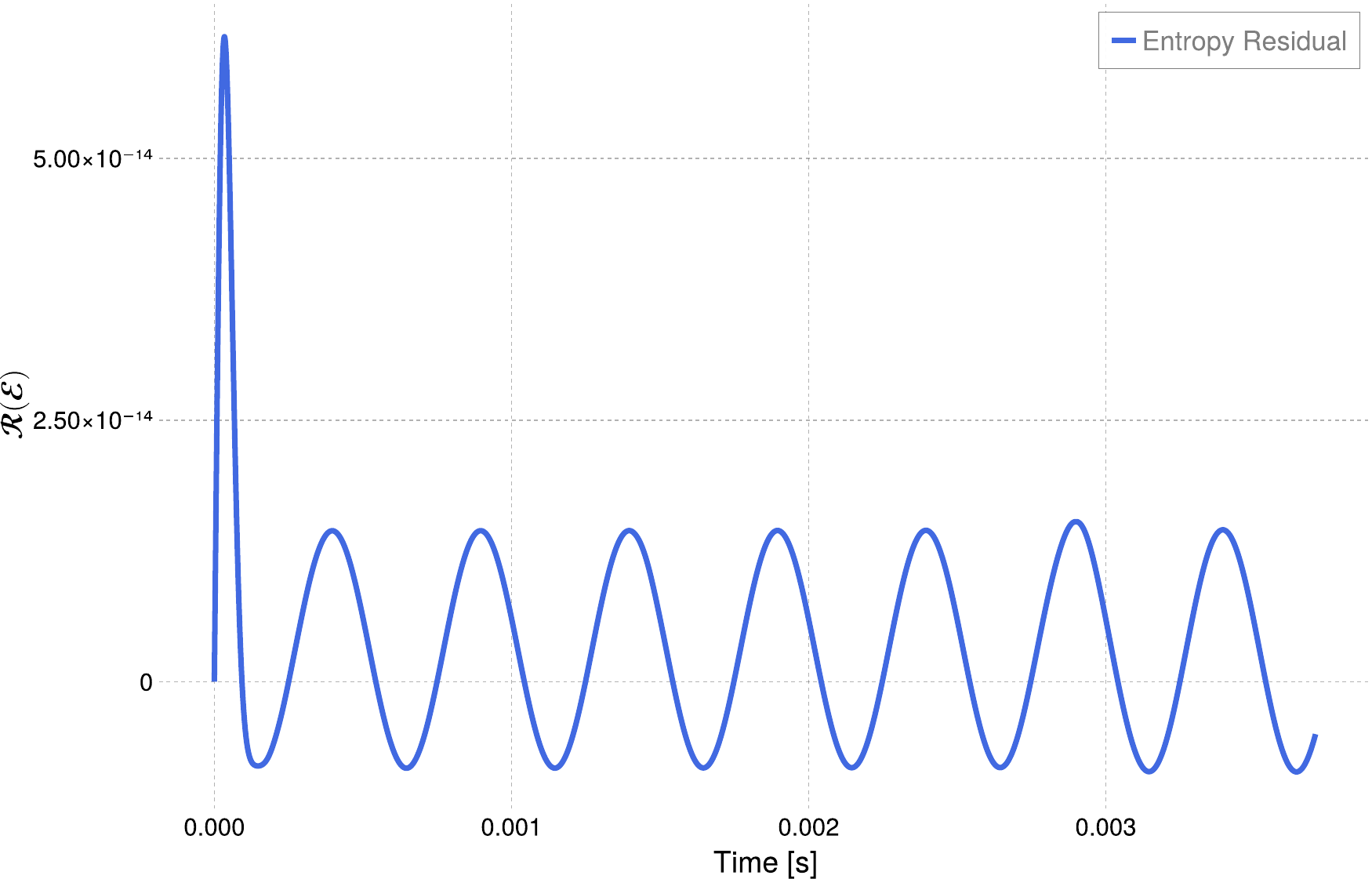}
\caption{Entropy residual versus time for the 3D test case, showing near-zero residual, indicating entropy conservation.}
\label{fig:3Dentropy-res}
\end{figure}

Fig.~\ref{fig:3Dentropy-res} shows the entropy residual for the 3D test case from Section~\ref{subsec:3D}. The entropy residual is close to zero, implying an entropy conservative scheme.

\end{appendices}

\bibliographystyle{unsrt}
\bibliography{references.bib}

\end{document}